\documentclass[11pt,letterpaper,reqno]{amsart}

\usepackage[T1]{fontenc}
\usepackage{amsmath,amssymb,amsfonts,amsthm,mathrsfs,bbm}
\usepackage{mathtools}
\usepackage{enumitem}
\usepackage{graphicx}
\usepackage{float}
\usepackage{booktabs}
\usepackage{doi}
\usepackage{bm}
\usepackage{tikz}
\usetikzlibrary{positioning, shapes.geometric, arrows.meta}
\usepackage{pgfplots}
\pgfplotsset{compat=1.18}
\usepackage{bookmark}
\usepackage{hyperref}
\hypersetup{pdfstartview={FitH}}

\newtheorem{thm}{Theorem}[section]
\newtheorem{lem}[thm]{Lemma}
\newtheorem{prop}[thm]{Proposition}
\newtheorem{cor}[thm]{Corollary}

\newtheorem{ques}[thm]{Question}
\newtheorem{conj}[thm]{Conjecture}

\theoremstyle{definition}

\newtheorem{rem}[thm]{Remark}

\numberwithin{equation}{section}
\allowdisplaybreaks[4]

\DeclareMathOperator{\Tr}{Tr}
\newcommand{\tr}{\Tr} 

\DeclareMathOperator{\rank}{rank}
\DeclareMathOperator{\supp}{supp}

\newcommand{\Id}{\mathrm{Id}}
\begin{document}
	
	\title[\(\ell_p\)--\(\ell_p\) and \(\ell_q\)--\(\ell_p\) Schatten 
	estimates]{Clarkson--McCarthy type inequalities, part I: \(\ell_p\)--\(\ell_p\) and \(\ell_q\)--\(\ell_p\) Schatten 
$p$-estimates}
	
	\author[T.~Zhang]{Teng Zhang}
	
	\address{School of Mathematics and Statistics, Xi'an Jiaotong University, Xi'an 710049, P. R. China}
	\email{teng.zhang@stu.xjtu.edu.cn}
	
	\subjclass[2020]{47A30, 46B20, 15A60}
	
	\keywords{Clarkson--McCarthy type inequalities; operator square-sum identities; complex interpolation}

	\begin{abstract}
		We characterize the matrices \(U=(u_{ij})\) for which the operator square-sum identity
		\[
		\sum_{i=1}^m\Big|\sum_{j=1}^n u_{ij}A_j\Big|^2=\sum_{j=1}^n|A_j|^2
		\]
		holds for all Schatten-class operators \(A_1,\ldots,A_n\); this happens exactly when \(U\) is an isometry.
		Using this characterization, we establish Clarkson--McCarthy type inequalities for several classes of operator families, including \(\ell_p\)--\(\ell_p\) estimates and mixed \(\ell_q\)--\(\ell_p\) estimates.
		We also obtain a multivariable extension of the Ball--Carlen--Lieb \(2\)-uniform convexity inequality and a weaker bound toward Audenaert's norm-compression conjecture.
	\end{abstract}

	\maketitle
	\tableofcontents
	
	\section{Introduction}
	For $1\le p<\infty$, let $S_p$ denote the Schatten--von Neumann class of compact operators
	$A:\ell_2\to\ell_2$ endowed with the norm
	\[
	\|A\|_p=\bigl(\tr (A^*A)^\frac{p}{2}\bigr)^\frac{1}{p}=\Bigl(\sum_{j} s_j(A)^p\Bigr)^\frac{1}{p},
	\]
	where $s_1(A)\ge s_2(A)\ge \cdots \ge 0$ denote the singular values of $A$, i.e., the eigenvalues of $|A|$
	arranged in nonincreasing order and counted with multiplicity.
	For $p=\infty$, we set
	\[
	\|A\|_\infty:=\sup_j s_j(A)=s_1(A).
	\]
	For $0<p<1$, we still define
	\[
	\|A\|_p=\Bigl(\sum_{j} s_j(A)^p\Bigr)^\frac{1}{p}<\infty,
	\]
	but this is only a quasi-norm (it does not satisfy the triangle inequality in general).
	
	In finite dimensions, we denote by $\mathbb{M}_{m,n}$ the space of $m\times n$ complex matrices; in particular,
	when $m=n=N$, we write $\mathbb{M}_N$. Since finite-rank operators are dense in $S_p$, in most proofs we may work
	in $\mathbb{M}_N$ and then pass to the limit. 
	
	\subsection{Clarkson--McCarthy inequalities}
	
	McCarthy \cite{Mc67} established the following two uniform convexity inequalities, namely, \emph{Clarkson--McCarthy inequalities}.
	\begin{thm}[McCarthy]\label{thm:McCarthy}
		Let $A,B\in S_p$. Then
		\begin{align}
			\|A+B\|_p^p+\|A-B\|_p^p &\le 2^{p-1}\big(\|A\|_p^p+\|B\|_p^p\big), \quad 2\le p<\infty, \label{e1}\\
			\|A+B\|_p^q+\|A-B\|_p^q &\le 2\big(\|A\|_p^p+\|B\|_p^p\big)^{\frac{q}{p}}, \qquad 1<p\le 2, \label{e2}
		\end{align}
		where $q$ denotes the conjugate exponent of $p$, i.e.\ $\frac1p+\frac1q=1$. Moreover, \eqref{e1} holds with the reverse inequality for $0<p\le 2$, and \eqref{e2} holds with the reverse inequality for $2\le p<\infty$.
	\end{thm}
	
	These two inequalities show that $S_p$ is $q$-uniformly convex for $1 < p \le 2$ and $p$-uniformly convex for $2\le p<\infty$ (compare with Clarkson's uniform convexity result for $L_p$ \cite{Cla36}). However, McCarthy was not the first to prove \eqref{e1}; this inequality already appeared in work of Dixmier \cite{Dix53}. Moreover, McCarthy's original proof of \eqref{e2} contains a gap, see the remark at the end of Fack--Kosaki \cite{FK86}, where the authors also provided a correct proof using complex interpolation.
	
	The Clarkson--McCarthy inequalities play an important role in analysis, operator theory, and mathematical physics, and have been studied extensively; see, for example, \cite{BCL94,BH88,BL20,FK86,HK08,Zha25}.
	
	Our motivation for this article stems from our desire to understand the connection between Theorem~\ref{thm:McCarthy} and the classical operator parallelogram law:
	\[
	\left| \frac{A+B}{\sqrt{2}}\right| ^2+\left| \frac{A-B}{\sqrt{2}}\right| ^2=|A|^2+|B|^2, \qquad A,B\in S_p.
	\]
	If we regard $(A, B)$ as an element of the product space $S_p\times S_p$, then we have the matrix representation
	$\begin{pmatrix}
		\frac{A+B}{\sqrt{2}}\\[4pt]
		\frac{A-B}{\sqrt{2}}
	\end{pmatrix}
	=
	\begin{pmatrix}
		\frac{\sqrt{2}}{2} & \frac{\sqrt{2}}{2}\\[4pt]
		\frac{\sqrt{2}}{2} & -\frac{\sqrt{2}}{2}
	\end{pmatrix}
	\begin{pmatrix}
		A\\
		B
	\end{pmatrix}$.
	This motivates us to seek analogues of Theorem~\ref{thm:McCarthy}. More precisely, suppose that there exists a matrix $U=(u_{ij})\in\mathbb{M}_2$ such that the following operator square-sum identity holds:
	\begin{align}\label{eq:operator_square_identity}
		\bigl|u_{11}A+u_{12}B\bigr|^2+\bigl|u_{21}A+u_{22}B\bigr|^2
		=
		|A|^2+|B|^2.
	\end{align}
	Can one then extend Theorem~\ref{thm:McCarthy} to this setting? We will see that if \eqref{eq:operator_square_identity} holds for every pair $(A,B)\in S_p\times S_p$, then $U$ must be a $2\times2$ unitary matrix (see Theorem~\ref{thm:identity_unitary}).
	
Ball, Carlen, and Lieb \cite{BCL94} established the following optimal
\(2\)-uniform convexity inequality, showing that for \(1< p\le 2\),
the space \(S_p\) is \(2\)-uniformly convex with sharp constant \(p-1\)
(compare with Hanner's uniform convexity result for \(L_p\) \cite{Han56}).
	
	\begin{thm}[Ball--Carlen--Lieb]\label{thm:Ball_Carlen_Lieb}
		Let $A,B\in S_p$. Then for $1\le p\le 2$,
		\begin{equation*}
			\|A+B\|_p^p+\|A-B\|_p^p\ge 2\left(  \|A\|_p^2+(p-1)\|B\|_p^2\right)^\frac{p}{2}.\end{equation*}
		For $2\le p<\infty$, the inequality is reversed.
	\end{thm}
	
	Theorem~\ref{thm:Ball_Carlen_Lieb} has many applications. For instance, it was generalized by King \cite{King03} to $2\times2$ block matrices. Motivated by problems in quantum information theory, King used his inequality to derive results on maximal $p$-norms of product channels. Moreover, $2$-uniform convexity is a cornerstone of Pisier's work \cite{PX03} and has further applications in mathematical physics (e.g., sharp fermion hypercontractivity~\cite{CL93}).

Operator Hanner's inequality, formerly known as the Ball--Carlen--Lieb conjecture, is the following noncommutative analogue of Hanner's classical inequality for \(L^p\)-spaces~\cite{Han56}.
	
	\begin{thm}[Ball--Carlen--Lieb--Hein\"avaara]\label{thm:BCLC}
		Let $A,B\in S_p$. For $1\le p\le 2$,
		\begin{equation}\label{e1.1}
			(\|A\|_p+\|B\|_p)^p+\bigl|\|A\|_p-\|B\|_p\bigr|^p\le \|A+B\|_p^p+\|A-B\|_p^p.
		\end{equation}
		For $2\le p<\infty$ the inequality is reversed.
	\end{thm}
	
	Ball, Carlen, and Lieb~\cite{BCL94} established \eqref{e1.1} for $1\le p\le 4/3$ and for $p\ge 4$, and they also proved it under additional positivity assumptions on $A\pm B$. The latter result was extended by King \cite{King03} to $2\times2$ block positive semidefinite matrices. A stronger variant proposed in \cite{CL06} was later disproved by Chayes \cite{Cha21}. Hein\"{a}vaara verified the case \(p=3\) in \cite{Hei23} and, more recently, established Theorem~\ref{thm:BCLC} in full in \cite{Hei25} by introducing the original and highly insightful notion of \emph{tracial joint spectral measures}. More precisely, Hein\"{a}vaara \cite{Hei25} proved that for fixed $A,B\in S_p$, the real linear span of $A$ and $B$ is isometric to a subspace of $L_p(\mu)$ for some positive measure $\mu$. However, his approach no longer works for more than two operators.
	
	\subsection{Clarkson--McCarthy type inequalities for several operators}
	Clarkson--McCarthy type inequalities with several operators have been intensely studied over the past few decades.
	
	In \cite{BK04}, Bhatia and Kittaneh established an $n$-operator version of Theorem~\ref{thm:McCarthy} with coefficients given by roots of unity.
	\begin{thm}[Bhatia--Kittaneh]\label{thm:Bhatia_Kittaneh}
		Let $A_1,\ldots,A_n\in S_p$. Then
		\begin{align}\label{bkp}
			\sum_{k=1}^{n}\left\| \sum_{j=1}^{n}\omega_{j-1}^{k-1}A_{j}\right\|_p^p&\le n^{p-1}\sum_{j=1}^{n}\left\| A_j\right\|_p^p,\quad   2\le p<\infty, \\
			\sum_{k=1}^{n}\left\| \sum_{j=1}^{n}\omega_{j-1}^{k-1}A_{j}\right\|_p^q&\le n\left( \sum_{j=1}^{n}\left\| A_j\right\|_p^p\right)^\frac{q}{p},\  1<p\le 2, \label{bkq}
		\end{align}
		where  $q$ denotes the conjugate exponent of $p$ and $\omega_0,\omega_{1},\ldots,\omega_{n-1}$ are the $n$-th roots of unity with $\omega_j=\exp(\tfrac{j2\pi i}{n})$. For $0<p\le 2$, \eqref{bkp} is reversed. For $2\le p<\infty$, \eqref{bkq} is reversed.
	\end{thm}
	
	In the next result, it is worth mentioning that Hirzallah and Kittaneh derived \eqref{hke1} in \cite{HK08}. Inequality \eqref{zhange1} (namely, \emph{Audenaert--Kittaneh's conjecture}) was conjectured by Audenaert and Kittaneh in \cite[Section~8]{AK17} and was recently proved by the author in \cite{Zha25}.
	\begin{thm}[Audenaert--Hirzallah--Kittaneh--Zhang]\label{thm:Audenaert_Hirazallah_Kittaneh_Zhang}
		Let $A_1,\ldots,A_n\in S_p$. Then
		\begin{align}
			\left\|\sum_{i=1}^nA_i\right\|_p^p+\sum_{1\le i<j\le n}\|A_i-A_j\|_p^p&\le n^{p-1}\sum_{j=1}^n\|A_j\|^p_p,\quad 2\le p<\infty,\label{hke1}\\
			\left\| \sum_{i=1}^nA_i\right\| _p^q+\sum_{1\le i<j\le n}\|A_i-A_j\|_p^q&\le n\left( \sum_{j=1}^n\|A_j\|^p_p\right)^\frac{q}{p},\ 1<p\le 2,\label{zhange1}
		\end{align}
		where $q$ stands for the conjugate exponent of $p$. For $0<p\le 2$, \eqref{hke1} is reversed. For $2\le p<\infty$, \eqref{zhange1} is reversed.
	\end{thm}
	
	Littlewood matrices are defined inductively by
	\[
	L_1=\begin{bmatrix} 1&1\\[2pt] 1&-1 \end{bmatrix},
	\qquad
	L_n=\begin{bmatrix}
		L_{n-1} & L_{n-1}\\
		L_{n-1} & -L_{n-1}
	\end{bmatrix}.
	\]
	They were first constructed by Sylvester in 1867. These matrices are often referred to as Walsh matrices and form a distinguished family of Hadamard matrices. Ke\v{c}ki\'{c} \cite[Corollary~4.5]{Ke19} established the following $2^n$-operator version of Theorem~\ref{thm:McCarthy}, with coefficients given by the entries of $L_n$.
	\begin{thm}[Ke\v{c}ki\'{c}]\label{thm:Keckic}
		Let $A_1,\ldots,A_{2^n}\in S_p$. Let  $\varepsilon_{ij}$ be the entries of the Littlewood matrix $L_n$. Then
		\begin{align}
			\sum_{i=1}^{2^n}\left\|\sum_{j=1}^{2^n}\varepsilon_{ij}A_j\right\|_p^p&\le 2^{n(p-1)}\sum\limits_{j=1}^{2^n}\|A_j\|_p^p,\quad  2\le p<\infty,\label{kep}\\
			\sum_{i=1}^{2^n}\left\|\sum_{j=1}^{2^n}\varepsilon_{ij}A_j\right\|_p^q&\le 2^{n}\left( \sum\limits_{j=1}^{2^n}\|A_j\|_p^p\right)^\frac{q}{p},\quad 1<p\le 2\label{keq},
		\end{align}
		where $q$ stands for the conjugate exponent of $p$. For $2\le p<\infty$, \eqref{keq} is reversed.
	\end{thm}
In fact, we will see that \eqref{kep} also reverses for $0<p\le 2$
(see Theorem~\ref{thm:schatten_p_norms}), a point that was not proved in
\cite{Ke19}.
	
	Let us now explain our findings. As noted above, these three families of Clarkson--McCarthy type inequalities for multiple operators all satisfy the following operator square-sum identity: for $A_1,\dots,A_n\in S_p$,
	\begin{align*}
		\sum_{k=1}^n \Big|\sum_{j=1}^n \omega_{j-1}^{k-1} A_j\Big|^2
		&= \sum_{j=1}^n \left|\sqrt{n}A_j\right|^2;
		\tag{Bhatia--Kittaneh}\\
		\Big|\sum_{i=1}^n A_i\Big|^2+\sum_{1\le i<j\le n}|A_i-A_j|^2
		&= \sum_{j=1}^n \left|\sqrt{n}A_j\right|^2;
		\tag{Hirzallah--Kittaneh}\\
		\sum_{i=1}^{2^n}\Big|\sum_{j=1}^{2^n}\varepsilon_{ij}A_j\Big|^2
		&= \sum_{j=1}^{2^n}\left|2^{\frac{n}{2}}A_j\right|^2. \tag{Ke\v{c}ki\'{c}}
	\end{align*}
This observation suggests that operator square-sum identities provide a natural framework for deriving multivariate Clarkson--McCarthy type inequalities.

In this paper, we first characterize the matrices \(U=(u_{ij})\in \mathbb{M}_{m,n}\) for which an operator square-sum identity of the form
\[
\sum_{i=1}^m\Big|\sum_{j=1}^n u_{ij}A_j\Big|^2=\sum_{j=1}^n|A_j|^2
\]
holds for all \(A_1,\ldots,A_n\in S_p\). We then use this characterization to obtain several families of Clarkson--McCarthy type inequalities for multiple operators. The basic mechanism is summarized in the following three theorems.
	
	\begin{thm}\label{thm:identity_unitary}
		For a fixed $m\times n$ matrix $U=(u_{ij})$,
		the operator square-sum identity
		\[
		\sum_{i=1}^m\Big| \sum_{j=1}^nu_{ij}A_j\Big| ^2=\sum_{j=1}^n|A_j|^2
		\]
		holds for all $A_1,\ldots,A_n\in S_p$
		if and only if $U$ is an isometry (i.e.\ $U^*U=I_n$).
	\end{thm}
	
	\begin{thm}\label{thm:schatten_p_norms}
		Let $B_i, A_j \in S_p$ for $1\le i\le m, 1\le j\le n$ satisfy
		\[
		\sum_{i=1}^m |B_i|^2 = \sum_{j=1}^n |A_j|^2,
		\]
		where $B_i$ are not necessarily linear combinations of $\{A_j\}$. Then for $2\le p<\infty$,
		\[
		m^{1 - \frac{p}{2}}\sum_{j=1}^n \|A_j\|_p^p \le \sum_{i=1}^m \|B_i\|_p^p \le n^{\frac{p}{2} - 1} \sum_{j=1}^n \|A_j\|_p^p.
		\]
		For $0 < p \le 2$, the inequalities are reversed.
	\end{thm}
	
	\begin{thm}\label{thm:schatten_p_q_norms}
		Let $B_i, A_j \in S_p$ for $1\le i\le m, 1\le j\le n$ satisfy $(B_1,\ldots,B_m)^T=U(A_1,\ldots,A_n)^T$ for some isometry $U=(u_{ij})$.  Then for $1<p\le 2$,
		\begin{align*}
			\sum_{i=1}^m \|B_i\|_p^q
			&\le
			(\max_{i,j}\left|u_{ij}\right|)^{\frac{q}{p}-1}
			\left(\sum_{j=1}^n \|A_j\|_p^p\right)^{\frac{q}{p}},\\
			\sum_{j=1}^n \|A_j\|_p^q
			&\le
			(\max_{i,j}\left|u_{ij}\right|)^{\frac{q}{p}-1}
			\left(	\sum_{i=1}^m \|B_i\|_p^p\right)^{\frac{q}{p}}.
		\end{align*}
		For $2\le p<\infty$, these two inequalities are reversed.
	\end{thm}
	
	\begin{rem}
		Theorems~\ref{thm:McCarthy}, \ref{thm:Bhatia_Kittaneh}, \ref{thm:Audenaert_Hirazallah_Kittaneh_Zhang} and~\ref{thm:Keckic} are special cases of Theorems~\ref{thm:schatten_p_norms} and~\ref{thm:schatten_p_q_norms}, obtained by choosing appropriate operator square-sum identities (as discussed above). Their complementary reverse-type inequalities are also given in Theorems~\ref{thm:schatten_p_norms} and~\ref{thm:schatten_p_q_norms}. For example,
		Zhang's 
		inequality~\eqref{zhange1} can be recovered from Theorem~\ref{thm:schatten_p_q_norms} by choosing the isometry
		\[
		U=\frac{1}{\sqrt{n}}\begin{pmatrix}
			1&1&1&\cdots&1&1\\
			1&-1&0&\cdots&0&0\\
			1&0&-1&\cdots&0&0\\
			\vdots&\vdots&\vdots&\ddots&\vdots&\vdots\\
			0&0&0&\cdots&1&-1
		\end{pmatrix}\in \mathbb{M}_{\frac{n(n-1)+2}{2}, n},
		\]
		whose rows encode the combinations $\sum_{i=1}^n A_i$ and $A_i-A_j$.
		Moreover, Theorem~\ref{thm:schatten_p_q_norms} yields the following sharp reverse-type companion to
		\eqref{zhange1}. 	\end{rem}
		\begin{cor}\label{cor:complementary-AKC}
			Let $A_1,\ldots,A_n\in S_p$. Then for $1<p\le 2$,
			\begin{equation}\label{eq:AKC-zhange}
				\sum_{j=1}^n\|A_j\|^q_p\le n^{-\frac{q}{p}}\left( \left\| \sum_{i=1}^nA_i\right\| _p^p+\sum_{1\le i<j\le n}\|A_i-A_j\|_p^p\right)^\frac{q}{p},
			\end{equation}
			where $q$ stands for the conjugate exponent of $p$. For $2\le p<\infty$, \eqref{eq:AKC-zhange} is reversed. Moreover, the constant $n^{-\frac{q}{p}}$ is optimal: equality holds for all $1<p<\infty$ by taking $A_1=\ldots=A_n \neq0$.
		\end{cor}

	As a corollary of Theorem~\ref{thm:schatten_p_norms}, we obtain a new proof of a classic result of Bhatia and Kittaneh~\cite{BK90} on block matrices, extending the admissible range to $p>0$. In the original paper \cite{BK90}, the result was proved only for $p\ge 1$, whereas our approach applies to all $p>0$.
	\begin{cor}[Bhatia--Kittaneh]\label{cor:BK90_newproof}
		Let $T=[T_{ij}]_{1\le i,j\le n}$ be an $n\times n$ block matrix whose blocks $T_{ij}$ belong to $\mathbb{M}_N$. Then for $2\le p<\infty$,
		\[
		n^{2-p}\|T\|_p^p \;\le\; \sum_{i,j=1}^n \|T_{ij}\|_p^p \;\le\; \|T\|_p^p.
		\]
		For $0< p\le 2$, both inequalities hold with the reverse sign.
	\end{cor}
	In February~2024,
	Carlen asked the author whether Theorem~\ref{thm:Ball_Carlen_Lieb} could be extended to
	the case of several operators. In this paper, by using a result
	of Ricard and Xu \cite{RX16}, we obtain the following generalization (in particular,
	when \(n=1\), it reduces to Theorem~\ref{thm:Ball_Carlen_Lieb}):
	\begin{thm}\label{thm:extension-BCL}
		Let $A,B_1,\ldots,B_n\in S_p$. Define $\varepsilon:=(\varepsilon_1, \ldots, \varepsilon_n)$. Then for $1<p\le 2$,
		\begin{equation}\label{eq:BCL-n-Sp-zhang}
			\frac{1}{2^n}\sum_{\varepsilon\in\{\pm1\}^n}
			\left\|A+\sum_{k=1}^n \varepsilon_k B_k\right\|_p^{p}
			\ge
			\left(\|A\|_p^{2}+(p-1)\sum_{k=1}^n \|B_k\|_p^{2}\right)^{\frac{p}{2}}.
		\end{equation}
		For $2\le p<\infty$, \eqref{eq:BCL-n-Sp-zhang} is reversed.
	\end{thm}

	As a generalization of Theorem~\ref{thm:BCLC}, Audenaert \cite{Aud08} (see also \cite[Section~5, Conjecture~3]{AK17}) proposed a conjecture concerning a norm compression inequality.
	\begin{conj}[Audenaert]\label{conj:Audenaert}
			Let $T $ be the $2\times n$ partitioned matrix
		\[
		T=\begin{pmatrix}
			A_1 & A_2 & \cdots & A_n\\
			B_1 & B_2 & \cdots & B_n
		\end{pmatrix},
		\]
		whose $(1,j)$-block is $A_j\in \mathbb M_N$ and whose $(2,j)$-block is $B_j\in \mathbb M_N$.
		Denoting by $C_p(T)$ its Schatten $p$-norm compression,
		\[
		C_p(T)=\begin{pmatrix}
			\|A_1\|_p & \|A_2\|_p & \cdots & \|A_n\|_p\\
			\|B_1\|_p & \|B_2\|_p & \cdots & \|B_n\|_p
		\end{pmatrix},
		\]
		we have
		\begin{align*}
			\|C_p(T)\|_p&\le	\|T\|_p , \qquad\qquad  1\le p\le 2, \\
			\|T\|_p &\le \|C_p(T)\|_p, \qquad 2\le p\le \infty.
		\end{align*}
		Moreover, 	the constants 1 in both inequalities are sharp.
	\end{conj}
To see that the constant \(1\) is sharp in both inequalities of Conjecture~\ref{conj:Audenaert}, let
	\[
	A_j=B_j=I_N,\qquad 1\le j\le n.
	\]
	Then
	\[
	T=\begin{pmatrix}
		I_N & I_N & \cdots & I_N\\
		I_N & I_N & \cdots & I_N
	\end{pmatrix}
	=
	\begin{pmatrix}
		1 & 1 & \cdots & 1\\
		1 & 1 & \cdots & 1
	\end{pmatrix}\otimes I_N,
	\]
	and
	\[
	C_p(T)=
	\begin{pmatrix}
		\|I_N\|_p & \|I_N\|_p & \cdots & \|I_N\|_p\\
		\|I_N\|_p & \|I_N\|_p & \cdots & \|I_N\|_p
	\end{pmatrix}
	=
	N^\frac{1}{p}
	\begin{pmatrix}
		1 & 1 & \cdots & 1\\
		1 & 1 & \cdots & 1
	\end{pmatrix}.
	\]
	Since the $2\times n$ matrix
	\[
	J=\begin{pmatrix}
		1 & 1 & \cdots & 1\\
		1 & 1 & \cdots & 1
	\end{pmatrix}
	\]
	has rank one and singular value $\sqrt{2n}$, we have
	\[
	\|J\|_p=\sqrt{2n}\qquad (1\le p\le \infty).
	\]
	Therefore
	\[
	\|T\|_p=\|J\otimes I_N\|_p=\|J\|_p\,\|I_N\|_p
	=\sqrt{2n}\,N^\frac{1}{p},
	\]
	while
	\[
	\|C_p(T)\|_p=N^\frac{1}{p}\|J\|_p
	=N^\frac{1}{p}\sqrt{2n}.
	\]
	Hence in this example
	\[
	\|C_p(T)\|_p=\|T\|_p.
	\]
	So the constant $1$ is attained, and therefore is sharp in both inequalities.
	
	It is easy to see that Conjecture~\ref{conj:Audenaert} implies Hanner's inequality (Theorem~\ref{thm:BCLC}) for matrices by taking
	$n=2$, $A_1=B_2=A$, and $A_2=B_1=B$. Unitarily conjugating $T$ by
	$\frac{1}{\sqrt{2}}
	\begin{pmatrix}
		I & I\\
		I & -I
	\end{pmatrix}$,
	and applying the corresponding norm compression with
	$
	\frac{1}{\sqrt{2}}
	\begin{pmatrix}
		1 & 1\\
		1 & -1
	\end{pmatrix}
	$
	yields Theorem~\ref{thm:BCLC} immediately.
	
	For $2\times 2$-partitioned block matrices $T$ (i.e., $n=2$), Conjecture~\ref{conj:Audenaert} is known to hold in two additional special cases: when $T$ is positive semidefinite \cite{King03}, and when all blocks of $T$ are diagonal matrices \cite{KN04}. Moreover, Conjecture~\ref{conj:Audenaert} has been established in the following situations~\cite{Aud08}: the norm compression of $T$ has rank~$1$; all blocks of $T$ have rank~$1$; all blocks $A_k$ in the first row are proportional, and likewise all blocks $B_k$ in the second row are proportional; and for general $2\times N$-partitioned $T$ whenever $p\ge 4$. 
	
We also remark that a \(3\times N\) analogue of Conjecture~\ref{conj:Audenaert} is false: Audenaert~\cite{Aud08} gave a \(3\times 3\) counterexample. Nevertheless, the next theorem shows that one still has a dimension-free estimate up to a rank-dependent factor.

\begin{thm}\label{thm:Audenaert-extension}
	Let $T$ be the $m\times n$ partitioned matrix
	\[
	T=\begin{pmatrix}
		A_{11} & A_{12} & \cdots & A_{1n}\\
		A_{21} & A_{22} & \cdots & A_{2n}\\
		\vdots&\vdots&\ddots&\vdots\\
		A_{m1} & A_{m2} & \cdots & A_{mn}
	\end{pmatrix},
	\]
	whose $(i,j)$-block is $A_{ij}\in \mathbb M_N$.
	Denoting by $C_p(T)$ its Schatten $p$-norm compression,
	\[
	C_p(T)=\begin{pmatrix}
		\|A_{11} \|_p & \|A_{12} \|_p & \cdots & \|A_{1n}\|_p\\
		\|A_{21} \|_p & \|A_{22} \|_p & \cdots & \|A_{2n}\|_p\\
		\vdots&\vdots&\ddots&\vdots\\
		\|A_{m1}\|_p & \|A_{m2}\|_p & \cdots & \|A_{mn}\|_p
	\end{pmatrix},
	\]
	we have
	\begin{align*}
		\|C_p(T)\|_p &\le r^{\frac1p-\frac12}\|T\|_p, \qquad\qquad 1\le p\le 2,\\
		\|T\|_p &\le r^{\frac12-\frac1p}\|C_p(T)\|_p, \qquad 2\le p\le \infty.
	\end{align*}
	where $r\le \min\{m,n\}$ is the rank of $C_p(T)$.
\end{thm}

Two immediate consequences of Theorem~\ref{thm:Audenaert-extension} are worth recording. First, in the case \(m=2\) we obtain the following weaker form of Conjecture~\ref{conj:Audenaert}.

\begin{cor}\label{cor:weaker_audenaert_conj}
	Let $T$ be the $2\times n$ partitioned matrix
	\[
	T=\begin{pmatrix}
		A_1 & A_2 & \cdots & A_n\\
		B_1 & B_2 & \cdots & B_n
	\end{pmatrix},
	\]
	whose $(1,j)$-block is $A_j\in \mathbb M_N$ and whose $(2,j)$-block is $B_j\in \mathbb M_N$.
	Denoting by $C_p(T)$ its Schatten $p$-norm compression,
	\[
	C_p(T)=\begin{pmatrix}
		\|A_1\|_p & \|A_2\|_p & \cdots & \|A_n\|_p\\
		\|B_1\|_p & \|B_2\|_p & \cdots & \|B_n\|_p
	\end{pmatrix},
	\]
	we have
	\begin{align*}
		\|C_p(T)\|_p &\le 2^{\frac1p-\frac12}\|T\|_p, \qquad\qquad 1\le p\le 2,\\
		\|T\|_p &\le 2^{\frac12-\frac1p}\|C_p(T)\|_p, \qquad 2\le p\le \infty.
	\end{align*}
\end{cor}

Second, when \(m=1\), the rank factor disappears and one recovers the one-row compression inequality with the sharp constant \(1\).

\begin{cor}\label{cor:one_row_compression}
	Let $T$ be the $1\times n$ partitioned matrix
	\[
	T=\begin{pmatrix}
		A_1 & A_2 & \cdots & A_n
	\end{pmatrix},
	\]
	whose $(1,j)$-block is $A_j\in \mathbb M_N$.
	Denoting by $C_p(T)$ its Schatten $p$-norm compression,
	\[
	C_p(T)=\begin{pmatrix}
		\|A_1\|_p & \|A_2\|_p & \cdots & \|A_n\|_p
	\end{pmatrix},
	\]
	we have
	\begin{align*}
		\|C_p(T)\|_p &\le \|T\|_p, \qquad\qquad 1\le p\le 2, \\
		\|T\|_p &\le \|C_p(T)\|_p, \qquad 2\le p\le \infty.
	\end{align*}
	Moreover, the constants \(1\) in both inequalities are sharp. Equality is attained by taking
$
	A_j=I_N
$
	for all $1\le j\le n$.
\end{cor}
\noindent\textbf{Organization of the paper.}
The paper is organized as follows. In Section~\ref{sec:examples}, we present some new examples of Clarkson--McCarthy type inequalities from the corresponding operator square-sum identities. In Section~\ref{sec:proof-iso}, we prove Theorems~\ref{thm:identity_unitary} and~\ref{thm:schatten_p_norms}, including the characterization of operator square-sum identities and the basic \(\ell_p\)--\(\ell_p\) Schatten estimates. Section~\ref{sec:proof-spq} is devoted to the proof of Theorem~\ref{thm:schatten_p_q_norms}, which yields the mixed \(\ell_q\)--\(\ell_p\) inequalities. In Section~\ref{sec:proof-bcl}, we prove Theorem~\ref{thm:extension-BCL}, establishing a multivariable extension of the Ball--Carlen--Lieb \(2\)-uniform convexity inequality. Finally, Section~\ref{sec:proof-aud} contains the proof of Theorem~\ref{thm:Audenaert-extension}; Corollaries~\ref{cor:weaker_audenaert_conj} and~\ref{cor:one_row_compression} then follow immediately.

\medskip
\noindent\textbf{Acknowledgements.}
The author is grateful to Professors Eric A.~Carlen and J.~C.~Bourin for their valuable comments and suggestions.
This work was supported by the China Scholarship Council, the Young Elite Scientists Sponsorship Program for PhD Students (China Association for Science and Technology), and the Fundamental Research Funds for the Central Universities at Xi'an Jiaotong University (Grant No.~xzy022024045).
\section{New examples}\label{sec:examples}

In this section, we record several new explicit isometries and the corresponding
Clarkson--McCarthy type inequalities furnished by
Theorems~\ref{thm:identity_unitary}, \ref{thm:schatten_p_norms}
and~\ref{thm:schatten_p_q_norms}.
\subsection{Euler's identity}
 Bourin and Lee \cite{BL26} asked for matrix analogues of Euler's quadrilateral identity.

\begin{ques}[{\cite[Question~4.10]{BL26}}]\label{ques:3}
	What are the matrix versions, if any, of Euler's quadrilateral identity
	\[
	\|x+y+z\|^{2}+\|x\|^{2}+\|y\|^{2}+\|z\|^{2}
	=
	\|x+y\|^{2}+\|y+z\|^{2}+\|z+x\|^{2}
	\]
	for three points $x,y,z\in\mathbb{C}^{n}$?
\end{ques}

A direct computation yields the following Euler-type operator identity.
Let $A,B,C\in S_p$. Then
	\begin{equation}\label{eq:Euler-op}
		|A+B+C|^{2}+|A|^{2}+|B|^{2}+|C|^{2}
		=
		|A+B|^{2}+|B+C|^{2}+|C+A|^{2}.
	\end{equation}

		Building on Euler's operator identity \eqref{eq:Euler-op}, we obtain the following sharp Clarkson--McCarthy type inequalities by Theorem~\ref{thm:schatten_p_norms}.

	\begin{prop}\label{prop:Euler-lp-lp}
		Let $A,B,C\in S_p$. Then for $2\le p<\infty$,
		\begin{equation*}
			\|A+B\|_p^p+\|B+C\|_p^p+\|C+A\|_p^p
			\le
			2^{p-2}\big(\|A+B+C\|_p^p+\|A\|_p^p+\|B\|_p^p+\|C\|_p^p\big).
		\end{equation*}
		For $0<p\le2$, the inequality  is reversed. Moreover, the constant $2^{p-2}$ is optimal: equality holds for all $p>0$ by taking $A=B=-C\neq0$.
	\end{prop}
	As a complement to Proposition~\ref{prop:Euler-lp-lp}, numerical evidence suggests the following sharp inequality. 
	For a weaker bound with coefficient $3^{\frac{p}{2}-1}$, it can be easily deduced by Theorem~\ref{thm:schatten_p_norms} based on \eqref{eq:Euler-op}.
	\begin{conj}\label{conj:sharp-inequality}
		Let $A,B,C\in S_p$.	Then for $2\le p<\infty$,
		\begin{equation*}
			\|A+B+C\|_p^p+\|A\|_p^p+\|B\|_p^p+\|C\|_p^p
			\le
			\frac{3^{p-1}+1}{2^p}\left(\|A+B\|_p^p+\|B+C\|_p^p+\|C+A\|_p^p\right).
		\end{equation*}
		For $0<p\le2$, the inequality  is reversed.  Moreover, the constant $\frac{3^{p-1}+1}{2^p}$ is optimal: equality holds for all $p>0$ by taking $A=B=C\neq0$.
	\end{conj}
	Classical \emph{Hlawka-type} inequalities compare norms (or their powers) of vector sums with those of their pairwise sums. 
	A standard form of Hlawka's inequality \cite{MPF93} asserts that if $V$ is an inner product space and $x,y,z\in V$, then
	\[
	\|x+y\|+\|y+z\|+\|z+x\|
	\ \le\ 
	\|x\|+\|y\|+\|z\|+\|x+y+z\|,
	\]
	where $\|\cdot\|$ denotes the norm induced by the inner product; equivalently, the inequality holds for any norm satisfying the parallelogram identity.
	Numerous refinements and variants replace $\|\cdot\|$ by powers $\|\cdot\|^{p}$ and/or adjust the constants.
	
Proposition~\ref{prop:Euler-lp-lp} and Conjecture~\ref{conj:sharp-inequality} can be viewed as a Schatten--$p$ (i.e.\ noncommutative $L_p$) analogue: it features a sharp $p$-dependent constant and exhibits a change in behavior at $p=2$, as is typical for inequalities governed by the uniform convexity/concavity of $L_p$ spaces.
	
	As a further refinement of Proposition~\ref{prop:Euler-lp-lp}, we establish the following sharp mixed $\ell_q$--$\ell_p$ Clarkson--McCarthy type inequalities.
	\begin{prop}\label{prop:Euler-lq-lp}
		Let $A,B,C\in S_p$.
		Then for $1<p\le 2$ ,
		\begin{align*}
			\|A+B+C\|_p^q+\|A\|_p^q+\|B\|_p^q+\|C\|_p^q	
			&\le
			2^{\,1-\frac{q}{p}}\,
			\Big(\|A+B\|_p^p+\|B+C\|_p^p+\|C+A\|_p^p\Big)^{\frac{q}{p}},\\
			\|A+B\|_p^q+\|B+C\|_p^q+\|C+A\|_p^q
			&\le
			2^{\,1-\frac{q}{p}}\,
			\Big(\|A+B+C\|_p^p+\|A\|_p^p+\|B\|_p^p+\|C\|_p^p	\Big)^{\frac{q}{p}},
		\end{align*}
		where $q$ is the conjugate exponent of $p$. For $2\le p<\infty$, these two inequalities are reversed. Moreover, the constant $2^{1-\frac{q}{p}}$ is optimal: equality holds for all $1< p<\infty$ by taking $A=B=-C\neq0$.
	\end{prop}
\begin{proof}
	Fix $1<p\le 2$ and let $q$ be the conjugate exponent of $p$.
	Set
	\[
	X_1:=A+B,\quad X_2:=B+C,\quad X_3:=C+A,
	\qquad
	Y_1:=A+B+C,\ Y_2:=A,\ Y_3:=B,\ Y_4:=C.
	\]
	Consider
	\[
	W:=\frac12
	\begin{pmatrix}
		1&1&1\\
		1&-1&1\\
		1&1&-1\\
		-1&1&1
	\end{pmatrix}\in\mathbb M_{4,3}.
	\]
	Then $W^*W=I_3$, so $W$ is an isometry, and $\max_{i,j}|w_{ij}|=1/2$.
	Moreover,
	\[
	\begin{pmatrix}
		Y_1\\Y_2\\Y_3\\Y_4
	\end{pmatrix}
	=
	W
	\begin{pmatrix}
		X_1\\X_2\\X_3
	\end{pmatrix},
	\]
	that is,
	\[
	Y_1=\frac{X_1+X_2+X_3}{2},\quad
	Y_2=\frac{X_1-X_2+X_3}{2},\quad
	Y_3=\frac{X_1+X_2-X_3}{2},\quad
	Y_4=\frac{-X_1+X_2+X_3}{2}.
	\]
	
	Applying Theorem~\ref{thm:schatten_p_q_norms} to $Y=WX$, we obtain
	\[
	\sum_{j=1}^4 \|Y_j\|_p^q
	\le
	\left(\frac12\right)^{\frac{q}{p}-1}
	\left(\sum_{i=1}^3 \|X_i\|_p^p\right)^{\frac{q}{p}},
	\]
	and
	\[
	\sum_{i=1}^3 \|X_i\|_p^q
	\le
	\left(\frac12\right)^{\frac{q}{p}-1}
	\left(\sum_{j=1}^4 \|Y_j\|_p^p\right)^{\frac{q}{p}}.
	\]
	This is exactly the desired pair of inequalities.
	
	For $2\le p<\infty$, the inequalities reverse by the reverse part of
	Theorem~\ref{thm:schatten_p_q_norms}.
	
	Finally, taking $A=B=-C\neq0$ yields equality in both inequalities, so the
	constant $2^{1-\frac{q}{p}}$ is optimal.
\end{proof}
	\begin{rem} 
		Setting $C=-B$ in Propositions~\ref{prop:Euler-lp-lp} and \ref{prop:Euler-lq-lp} yields the classical Clarkson--McCarthy inequalities (Theorem~\ref{thm:McCarthy}).
	\end{rem}

\subsection{A mean-variance identity}
Let $x_1, x_2, \dots, x_n$ be real numbers associated with probabilities $p_1, p_2, \dots, p_n$, where $p_i \ge 0$ for all $i$ and $\sum_{i=1}^n p_i = 1$. Define the mean by
\[
\mathsf{m}= \sum_{i=1}^n p_i x_i .
\]
Then the variance is given by
\[
\sum_{j=1}^n p_j (x_j-\mathsf{m})^2.
\]
Expanding this expression yields
\[
\sum_{j=1}^n p_j (x_j-\mathsf{m})^2
=
\sum_{j=1}^n p_j x_j^2 -\mathsf{m}^2,
\]
and hence
\begin{equation}\label{eq:mean-variance}
	\sum_{j=1}^n p_j x_j^2
	=
	\sum_{j=1}^n p_j (x_j-\mathsf{m})^2 + \mathsf{m}^2.	
\end{equation}
This is exactly the discrete probabilistic form of the mean-variance identity
\[
\mathbb{E}[X^2] = \operatorname{Var}(X) + \bigl(\mathbb{E}[X]\bigr)^2,
\]
where the random variable $X$ takes the value $x_i$ with probability $p_i$.

The operator analogue of \eqref{eq:mean-variance} is given by
\begin{equation}\label{eq:operator-mean-variance-prop}
	\left| \sum_{i=1}^n p_i A_i\right| ^2+	\sum_{j=1}^n p_j \left| A_j-\sum_{i=1}^n p_i A_i\right| ^2
	=\sum_{j=1}^n p_j |A_j|^2.
\end{equation}
We now derive Clarkson-type inequalities from the identity \eqref{eq:operator-mean-variance-prop}.

\begin{prop}\label{prop:operator_mean_variance_estimate}
	Let $A_1,\ldots,A_n\in S_p$ and let $p_1,\ldots,p_n\ge 0$ satisfy
$
	\sum_{j=1}^n p_j=1.
$
	Then
	\begin{align}
		\Bigl\|\sum_{i=1}^n p_iA_i\Bigr\|_p^p
		+\sum_{j=1}^n p_j^{\frac{p}{2}}\Bigl\|A_j-\sum_{i=1}^n p_iA_i\Bigr\|_p^p
		&\le
		n^{\frac{p}{2}-1}\sum_{j=1}^n p_j^{\frac{p}{2}}\|A_j\|_p^p,
		\qquad	\qquad	\qquad2\le p<\infty, \label{eq:omv1}\\
		\Bigl\|\sum_{i=1}^n p_iA_i\Bigr\|_p^q
		+\sum_{j=1}^n p_j^{\frac{q}{2}}\Bigl\|A_j-\sum_{i=1}^n p_iA_i\Bigr\|_p^q
		&\le
		\gamma_n^{\frac{q}{p}-1}
		\left(\sum_{j=1}^n p_j^{\frac{p}{2}}\|A_j\|_p^p\right)^{\frac{q}{p}}, 1<p\le 2, \label{eq:omv2}
	\end{align}
	where  $q$ denotes the conjugate exponent of $p$ and $	\gamma_n
	:=
	\max\limits\left\{
	\max\limits_{1\le j\le n}\sqrt{p_j},
	\max\limits_{1\le i\le n}(1-p_i)
	\right\}$. Moreover, \eqref{eq:omv1} holds with the reverse inequality for $0<p\le 2$, and \eqref{eq:omv2} holds with the reverse inequality for $2\le p<\infty$.
\end{prop}

\begin{proof}
Clearly, \eqref{eq:operator-mean-variance-prop} can be rewritten as
	\[
	\left|\sum_{i=1}^n p_iA_i\right|^2
	+\sum_{j=1}^n\left|\sqrt{ p_j}\left( A_j-\sum_{i=1}^n p_iA_i\right) \right|^2
	=
	\sum_{j=1}^n|\sqrt{ p_j}A_j|^2.
	\]
	Hence Theorem~\ref{thm:schatten_p_norms}
	yields
	\[
	\Bigl\|\sum_{i=1}^n p_iA_i\Bigr\|_p^p
	+\sum_{j=1}^n p_j^{\frac{p}{2}}\Bigl\|A_j-\sum_{i=1}^n p_iA_i\Bigr\|_p^p
	\le
	n^{\frac{p}{2}-1}\sum_{j=1}^n p_j^{\frac{p}{2}}\|A_j\|_p^p,
	\qquad 2\le p<\infty,
	\]
	with the reverse inequality for \(0<p\le 2\). This proves \eqref{eq:omv1}.
	
	For \eqref{eq:omv2}, consider the \((n+1)\times n\) matrix
	\[
	U=
	\begin{pmatrix}
		\sqrt{p_1} & \sqrt{p_2} & \cdots & \sqrt{p_n}\\
		1-p_1 & -\sqrt{p_1p_2} & \cdots & -\sqrt{p_1p_n}\\
		-\sqrt{p_2p_1} & 1-p_2 & \cdots & -\sqrt{p_2p_n}\\
		\vdots & \vdots & \ddots & \vdots\\
		-\sqrt{p_np_1} & -\sqrt{p_np_2} & \cdots & 1-p_n
	\end{pmatrix}.
	\]
	Then
	\[
	U^*U=I_n,
	\]
	so \(U\) is an isometry, and
	\[
	\begin{pmatrix}
		\sum_{i=1}^n p_iA_i\\[1mm]
		\sqrt{p_1}\left(A_1-\sum_{i=1}^n p_iA_i\right)\\
		\vdots\\
		\sqrt{p_n}\left(A_n-\sum_{i=1}^n p_iA_i\right)
	\end{pmatrix}
	=
	U
	\begin{pmatrix}
		\sqrt{p_1}A_1\\
		\vdots\\
		\sqrt{p_n}A_n
	\end{pmatrix}.
	\]
	Therefore, Theorem~\ref{thm:schatten_p_q_norms} gives
	\begin{align*}
		\Bigl\|\sum_{i=1}^n p_iA_i\Bigr\|_p^q
		+\sum_{j=1}^n p_j^{\frac{q}{2}}\Bigl\|A_j-\sum_{i=1}^n p_iA_i\Bigr\|_p^q
		\le
	(\max_{i,j}|u_{ij}|)^{\frac{q}{p}-1}
		\left(\sum_{j=1}^n p_j^{\frac{p}{2}}\|A_j\|_p^p\right)^{\frac{q}{p}},
		\ 1<p\le 2,
	\end{align*}
	with the reverse inequality for \(2\le p<\infty\). Since
	\[
	u_{1j}=\sqrt{p_j},
	\qquad
	u_{ij}=
	\begin{cases}
		1-p_{i-1}, & i=j+1,\\
		-\sqrt{p_{i-1}p_j}, & i\ne j+1,
	\end{cases}
	\]
	we obtain
	\[
	\max_{i,j}|u_{ij}|
	=
	\max\left\{
	\max_{1\le j\le n}\sqrt{p_j},
	\max_{1\le i\le n}(1-p_i)
	\right\}.
	\]
	This proves \eqref{eq:omv2}.
\end{proof}

\begin{rem}
	In particular, when $p_1=\cdots=p_n=\frac1n$, \eqref{eq:omv1} becomes
	\[
	\Bigl\|\frac1n\sum_{i=1}^n A_i\Bigr\|_p^p
	+n^{-\frac{p}{2}}\sum_{j=1}^n \Bigl\|A_j-\frac1n\sum_{i=1}^n A_i\Bigr\|_p^p
	\le
	\frac1n\sum_{j=1}^n \|A_j\|_p^p,
	\qquad 2\le p<\infty,
	\]
	with the reverse inequality for $0<p\le 2$. Moreover, the constant $\frac1n$ on the right-hand side is sharp. Indeed, if $A_1=\cdots=A_n=A\neq 0$, then
 equality is attained.
	Therefore, the constant $\frac1n$ cannot be improved.
\end{rem}

\subsection{All sign combinations}
For each $\varepsilon=(\varepsilon_1,\ldots,\varepsilon_n)\in\{\pm1\}^n$, define
\[
S_\varepsilon:=\sum_{k=1}^n \varepsilon_k A_k.
\]
By summing over all sign combinations, the mixed terms cancel, and one obtains the identity
\[
\sum_{\varepsilon\in\{\pm1\}^n}|S_\varepsilon|^2
=
\sum_{k=1}^n |2^{\frac{n}{2}}A_k|^2.
\]
This identity is induced by an explicit isometry, and hence gives rise to the following Clarkson--McCarthy type inequalities.

\begin{prop}\label{prop:all-signs} Let \(A_1,\ldots,A_n\in S_p\).
	For each $\varepsilon=(\varepsilon_1,\ldots,\varepsilon_n)\in\{\pm1\}^n$,  set
$
	S_\varepsilon:=\sum_{k=1}^n \varepsilon_k A_k.
$
	Then for $2\le p<\infty$,
	\[
	2^n\sum_{k=1}^n \|A_k\|_p^p
	\le
	\sum_{\varepsilon\in\{\pm1\}^n}\|S_\varepsilon\|_p^p
	\le
	2^{n\frac{p}{2}}n^{\frac{p}{2}-1}\sum_{k=1}^n \|A_k\|_p^p,
	\]
	with reversed inequalities for $0<p\le 2$.
	If $1<p\le 2$ and $q$ is conjugate to $p$, then
\begin{align*}
		\sum_{\varepsilon\in\{\pm1\}^n}\|S_\varepsilon\|_p^q
	&\le
	2^n\left(\sum_{k=1}^n \|A_k\|_p^p\right)^{\frac{q}{p}},\\
	\sum_{k=1}^n \|A_k\|_p^q
	&\le
	2^{-n\frac{q}{p}}
	\left(
	\sum_{\varepsilon\in\{\pm1\}^n}\|S_\varepsilon\|_p^p
	\right)^{\frac{q}{p}},
\end{align*}
	with reversed inequalities for $2\le p<\infty$.
\end{prop}

\begin{proof}
	Let $U\in\mathbb M_{2^n,n}$ be defined by
	\[
	u_{\varepsilon,k}=2^{-\frac{n}{2}}\varepsilon_k,
	\qquad
	\varepsilon\in\{\pm1\}^n,\ \ 1\le k\le n.
	\]
	Then
	\[
	(U^*U)_{jk}
	=
	2^{-n}\sum_{\varepsilon\in\{\pm1\}^n}\varepsilon_j\varepsilon_k
	=
	\delta_{jk},
	\]
	so $U$ is an isometry. Therefore, by Theorem~\ref{thm:identity_unitary},
	\[
	\sum_{\varepsilon\in\{\pm1\}^n}|S_\varepsilon|^2
	=
	\sum_{k=1}^n |2^{\frac{n}{2}}A_k|^2.
	\]
	Applying Theorem~\ref{thm:schatten_p_norms} to the families
	$\{S_\varepsilon\}_{\varepsilon\in\{\pm1\}^n}$ and
	$\{2^{\frac{n}{2}}A_k\}_{k=1}^n$, we obtain, for $2\le p<\infty$,
	\[
	2^n\sum_{k=1}^n \|A_k\|_p^p
	\le
	\sum_{\varepsilon\in\{\pm1\}^n}\|S_\varepsilon\|_p^p
	\le
	2^{n\frac p2}n^{\frac p2-1}\sum_{k=1}^n \|A_k\|_p^p,
	\]
	with reversed inequalities for $0<p\le 2$.
	
	Moreover, since
	\[
	\max_{\varepsilon,k}|u_{\varepsilon,k}|=2^{-\frac{n}{2}},
	\]
	Theorem~\ref{thm:schatten_p_q_norms} yields, for $1<p\le 2$,
	\[
	\sum_{\varepsilon\in\{\pm1\}^n}\|S_\varepsilon\|_p^q
	\le
	2^n\left(\sum_{k=1}^n \|A_k\|_p^p\right)^{\frac{q}{p}},
	\]
	and
	\[
	2^{n\frac{q}{2}}\sum_{k=1}^n \|A_k\|_p^q
	\le
	2^{-\frac n2(\frac qp-1)}
	\left(
	\sum_{\varepsilon\in\{\pm1\}^n}\|S_\varepsilon\|_p^p
	\right)^{\frac{q}{p}},
	\]
	which is equivalent to
	\[
	\sum_{k=1}^n \|A_k\|_p^q
	\le
	2^{-n\frac{q}{p}}
	\left(
	\sum_{\varepsilon\in\{\pm1\}^n}\|S_\varepsilon\|_p^p
	\right)^{\frac{q}{p}}.
	\]
	For $2\le p<\infty$, the two inequalities are reversed.
\end{proof}

\subsection{A star-type identity}
A direct computation shows that
\[
\left|\sum_{i=1}^n A_i\right|^2
+
\sum_{j=1}^n\left|\sum_{i=1}^n A_i-\frac{n+1}{2}A_j\right|^2
=
\sum_{j=1}^n\left|\frac{n+1}{2}A_j\right|^2.
\]
This identity is again generated by an explicit isometry, and therefore leads to the following Clarkson--McCarthy type inequalities.

\begin{prop}\label{prop:star}
Let $A_1,\ldots,A_n\in S_p$.
	Then for $2\le p<\infty$,
	\[
	(n+1)^{1-\frac{p}{2}}\sum_{j=1}^n \|A_j\|_p^p
	\le
	\left\|\frac{2}{n+1}\sum_{i=1}^n A_i\right\|_p^p
	+
	\sum_{j=1}^n \left\|\frac{2}{n+1}\sum_{i=1}^n A_i-A_j\right\|_p^p
	\le
	n^{\frac{p}{2}-1}\sum_{j=1}^n \|A_j\|_p^p,
	\]
	with reversed inequalities for $0<p\le 2$.
	If $1<p\le 2$ and $q$ is conjugate to $p$, then
	\[
	\left\|\frac{2}{n+1}\sum_{i=1}^n A_i\right\|_p^q
	+
	\sum_{j=1}^n \left\|\frac{2}{n+1}\sum_{i=1}^n A_i-A_j\right\|_p^q
	\le
	\beta_n^{\frac{q}{p}-1}
	\left(\sum_{j=1}^n \|A_j\|_p^p\right)^{\frac{q}{p}},
	\]
	and
	\[
	\sum_{j=1}^n \|A_j\|_p^q
	\le
	\beta_n^{\frac{q}{p}-1}
	\left(
	\left\|\frac{2}{n+1}\sum_{i=1}^n A_i\right\|_p^p
	+
	\sum_{j=1}^n \left\|\frac{2}{n+1}\sum_{i=1}^n A_i-A_j\right\|_p^p
	\right)^{\frac{q}{p}},
	\]
	with reversed inequalities for $2\le p<\infty$, where $	\beta_n:=\max\left\{\frac{2}{n+1},\frac{n-1}{n+1}\right\}.$
\end{prop}

\begin{proof}
	Let $\mathbf 1=(1,\ldots,1)^T\in\mathbb C^n$, let $e_1,\ldots,e_n$ be the standard basis,
	and define
	\[
	v_0:=\frac{2}{n+1}\mathbf 1,
	\qquad
	v_j:=\frac{2}{n+1}\mathbf 1-e_j
	\quad (1\le j\le n).
	\]
	Let $U\in\mathbb M_{n+1,n}$ be the matrix whose rows are
	$v_0^T,v_1^T,\ldots,v_n^T$.
	A direct computation shows that
	\[
	U^*U=I_n.
	\]
Hence \(U\) is an isometry. Set
\[
S:=\sum_{i=1}^n A_i.
\]
By Theorem~\ref{thm:identity_unitary}, we have
\[
\left|\frac{2}{n+1}S\right|^2
+
\sum_{j=1}^n\left|\frac{2}{n+1}S-A_j\right|^2
=
\sum_{j=1}^n |A_j|^2.
\]
Multiplying both sides by \(\bigl(\frac{n+1}{2}\bigr)^2\) gives the identity displayed just before Proposition~\ref{prop:star}.
	
	Applying Theorem~\ref{thm:schatten_p_norms} to the families
	\[
	B_0:=\frac{2}{n+1}S,
	\qquad
	B_j:=\frac{2}{n+1}S-A_j\quad (1\le j\le n),
	\]
	and $\{A_j\}_{j=1}^n$, we obtain the $\ell_p\to\ell_p $  inequalities.
	
	Since the entries of $U$ are either $\frac{2}{n+1}$ or $-\frac{n-1}{n+1}$, we have
	\[
	\max_{i,j}|u_{ij}|
	=
	\beta_n
	=
	\max\left\{\frac{2}{n+1},\frac{n-1}{n+1}\right\}.
	\]
	Hence Theorem~\ref{thm:schatten_p_q_norms} yields the  $\ell_q\to\ell_p $  inequalities.
\end{proof}
	\section{Proof of Theorems~\ref{thm:identity_unitary} and \ref{thm:schatten_p_norms}}\label{sec:proof-iso}
	
	\begin{proof}[Proof of Theorem~\ref{thm:identity_unitary}]
		Let $U=(u_{ij})\in \mathbb M_{m,n}$ be fixed and suppose that
		\begin{equation}\label{eq:squaresum}
			\sum_{i=1}^m\left|\sum_{j=1}^n u_{ij}A_j\right|^2=\sum_{j=1}^n|A_j|^2
			\qquad\text{for all }A_1,\dots,A_n\in S_p.
		\end{equation}
		
		\medskip
		\noindent\textbf{Necessity.}
		Fix a nonzero finite-rank projection $P$ (e.g.\ rank one). For arbitrary scalars $\alpha_1,\dots,\alpha_n\in\mathbb C$, set
		$A_j=\alpha_j P$. Then $|A_j|^2=|\alpha_j|^2P$ and
		\[
		\left|\sum_{j=1}^n u_{ij}A_j\right|^2
		=
		\left|\sum_{j=1}^n u_{ij}\alpha_j\right|^2 P.
		\]
		Substituting into \eqref{eq:squaresum}, we obtain
		\[
		\sum_{i=1}^m\left|\sum_{j=1}^n u_{ij}\alpha_j\right|^2 P
		=
		\left(\sum_{j=1}^n|\alpha_j|^2\right)P.
		\]
		Since $P\neq 0$, equality of these scalar multiples of $P$ implies, for all
		$\alpha=(\alpha_j)\in\mathbb C^n$,
		\[
		\sum_{i=1}^m\left|\sum_{j=1}^n u_{ij}\alpha_j\right|^2
		=
		\sum_{j=1}^n|\alpha_j|^2,
		\]
		i.e.\ $\|U\alpha\|_2=\|\alpha\|_2$. Hence $U^*U=I_n$ and $U$ is an isometry.
		
		\medskip
		\noindent\textbf{Sufficiency.}
		If $U^*U=I_n$, then expanding as in the standard computation gives
		\[
		\sum_{i=1}^m\left|\sum_{j=1}^n u_{ij}A_j\right|^2
		=
		\sum_{j,k=1}^n (U^*U)_{jk}A_j^*A_k=\sum_{j=1}^n A_j^*A_j=\sum_{j=1}^n|A_j|^2,
		\]
		which is \eqref{eq:squaresum}.
	\end{proof}
In the proof of Theorem~\ref{thm:schatten_p_norms}, we shall use the following classical trace inequality due to Rotfel'd for a finite sum of positive operators; it is an immediate consequence of \cite[Theorem~1]{Rot69}.

\begin{lem}[\cite{Rot69}]\label{lem:rotfeld}
	Let $f:[0,\infty)\to[0,\infty)$ be continuous with $f(0)=0$. For any finite family of positive operators $\{X_k\}_{k=1}^n$ the following hold:
	\begin{enumerate}
		\item If $f$ is \emph{concave}, then
		\[
		\Tr f\Big(\sum_{k=1}^n X_k\Big)\ \le\ \sum_{k=1}^n \Tr f(X_k).
		\]
		\item If $f$ is \emph{convex}, then
		\[
		\Tr f\Big(\sum_{k=1}^n X_k\Big)\ \ge\ \sum_{k=1}^n \Tr f(X_k).
		\]
	\end{enumerate}
\end{lem}
Next, we prove Theorem~\ref{thm:schatten_p_norms}.
\begin{proof}[Proof of Theorem~\ref{thm:schatten_p_norms}]
	Set
	\[
	P_i:=|B_i|^2,\qquad Q_j:=|A_j|^2,
	\]
	so $P_i,Q_j\ge 0$ and
\begin{equation}\label{eq:sum-Pi-Qj-equal-S}
	\sum_{i=1}^m P_i=\sum_{j=1}^n Q_j=:S\ge 0.
\end{equation}
	Let $r:=\frac{p}{2}$. For any $X\in S_p$, we have
	\[
	\|X\|_p^p=\Tr(|X|^p)=\Tr\big((|X|^2)^{\frac{p}{2}}\big)=\Tr\big((|X|^2)^r\big),
	\]
	hence
\begin{equation}\label{eq:sum-Bi-Ai}
		\sum_{i=1}^m \|B_i\|_p^p=\sum_{i=1}^m \Tr(P_i^r),
	\qquad
	\sum_{j=1}^n \|A_j\|_p^p=\sum_{j=1}^n \Tr(Q_j^r).
\end{equation}

	\medskip
	\noindent\textbf{Step 1: upper bound for $2\le p<\infty$ (i.e.\ $r\ge 1$).}
	Since $t\mapsto t^r$ is convex on $[0,\infty)$ and $t^r(0)=0$, Lemma~\ref{lem:rotfeld} (convex case) and \eqref{eq:sum-Pi-Qj-equal-S}
	yields, for positive operators,
\begin{equation}\label{eq:Pi-less-S}
	\sum_{i=1}^m \Tr(P_i^r)\ \le\ \Tr\Big(\sum_{i=1}^m P_i\Big)^r
=\Tr(S^r).
\end{equation}
	On the other hand, using the triangle inequality in $S_r$ and the scalar estimate
	\[
	\Big(\sum_{j=1}^n a_j\Big)^r\le n^{r-1}\sum_{j=1}^n a_j^r
	\qquad(a_j\ge 0,\ r\ge 1),
	\]
by \eqref{eq:sum-Pi-Qj-equal-S}, we obtain
\begin{equation}\label{eq:S-lee-Qj}
		\Tr(S^r)=\Big\|\sum_{j=1}^n Q_j\Big\|_r^r
	\le \Big(\sum_{j=1}^n \|Q_j\|_r\Big)^r
	\le n^{r-1}\sum_{j=1}^n \|Q_j\|_r^r
	= n^{r-1}\sum_{j=1}^n \Tr(Q_j^r).
\end{equation}
	Combining \eqref{eq:Pi-less-S} and \eqref{eq:S-lee-Qj} gives
	\[
	\sum_{i=1}^m \Tr(P_i^r)\ \le\ n^{r-1}\sum_{j=1}^n \Tr(Q_j^r),
	\]
	by \eqref{eq:sum-Bi-Ai}, which is exactly
	\[
	\sum_{i=1}^m \|B_i\|_p^p \le n^{\frac p2-1}\sum_{j=1}^n \|A_j\|_p^p.
	\]
	
	\medskip
	\noindent\textbf{Step 2: lower bound for $2\le p<\infty$.}
	Applying Lemma~\ref{lem:rotfeld} (convex case) to $\{Q_j\}$ and \eqref{eq:sum-Pi-Qj-equal-S} gives
\begin{equation}\label{eq:lower-1}
	\sum_{j=1}^n \Tr(Q_j^r)\ \le\ \Tr\Big(\sum_{j=1}^n Q_j\Big)^r
	=\Tr(S^r).
\end{equation}
Similar to \eqref{eq:S-lee-Qj},	applying the previous $S_r$ estimate to \eqref{eq:sum-Pi-Qj-equal-S} yields
\begin{equation}\label{eq:lower-2}
		\Tr(S^r)=\Big\|\sum_{i=1}^m P_i\Big\|_r^r
	\le m^{r-1}\sum_{i=1}^m \|P_i\|_r^r
	= m^{r-1}\sum_{i=1}^m \Tr(P_i^r).
\end{equation}
Combining \eqref{eq:lower-1} and \eqref{eq:lower-2} gives
	\[
	\sum_{j=1}^n \Tr(Q_j^r)
	\le m^{r-1}\sum_{i=1}^m \Tr(P_i^r),
	\]
	i.e.
	\[
	m^{1-r}\sum_{j=1}^n \Tr(Q_j^r)\ \le\ \sum_{i=1}^m \Tr(P_i^r),
	\]
	equivalently, by \eqref{eq:sum-Bi-Ai}
	\[
	m^{1-\frac p2}\sum_{j=1}^n \|A_j\|_p^p
	\le
	\sum_{i=1}^m \|B_i\|_p^p.
	\]
	This completes the proof for $2\le p<\infty$.
	
\medskip
\noindent\textbf{Step 3: the case $0<p\le 2$.}
Here $r=\frac{p}{2}\in(0,1]$, so $t\mapsto t^r$ is concave on $[0,\infty)$.
By Lemma~\ref{lem:rotfeld} and \eqref{eq:sum-Pi-Qj-equal-S},
\begin{equation}\label{eq:inverse1}
	\sum_{i=1}^m \Tr(P_i^r)\ge \Tr\Big(\sum_{i=1}^m P_i\Big)^r
	=\Tr(S^r),
	\qquad
	\sum_{j=1}^n \Tr(Q_j^r)\ge \Tr\Big(\sum_{j=1}^n Q_j\Big)^r
	=\Tr(S^r).
\end{equation}

To recover the quantitative constants, we use the fact that $t\mapsto t^r$ is
operator concave on $[0,\infty)$ for $0<r\le 1$. Hence
\[
\Bigl(\frac1n\sum_{j=1}^n Q_j\Bigr)^r
\ge
\frac1n\sum_{j=1}^n Q_j^r,
\]
and therefore
\begin{equation}\label{eq:inverse2}
	\Tr(S^r)
	=
	n^r\Tr\Bigl(\frac1n\sum_{j=1}^n Q_j\Bigr)^r
	\ge
	n^{r-1}\sum_{j=1}^n \Tr(Q_j^r).
\end{equation}
Similarly,
\begin{equation}\label{eq:inverse3}
\Tr(S^r)
=
m^r\Tr\Bigl(\frac1m\sum_{i=1}^m P_i\Bigr)^r
\ge
m^{r-1}\sum_{i=1}^m \Tr(P_i^r).
\end{equation}
Combining \eqref{eq:inverse1}, \eqref{eq:inverse2} and \eqref{eq:inverse3}, we obtain
\[
n^{r-1}\sum_{j=1}^n \Tr(Q_j^r)
\le
\Tr(S^r)
\le
\sum_{i=1}^m \Tr(P_i^r),
\]
and
\[
m^{r-1}\sum_{i=1}^m \Tr(P_i^r)
\le
\Tr(S^r)
\le
\sum_{j=1}^n \Tr(Q_j^r).
\]
Equivalently, by \eqref{eq:sum-Bi-Ai}
\[
\sum_{i=1}^m \|B_i\|_p^p
\ge
n^{\frac p2-1}\sum_{j=1}^n \|A_j\|_p^p,
\qquad
\sum_{i=1}^m \|B_i\|_p^p
\le
m^{1-\frac p2}\sum_{j=1}^n \|A_j\|_p^p.
\]
This is exactly the reversed form of the theorem for $0<p\le 2$.
\end{proof}
As a direct application of Theorem~\ref{thm:schatten_p_norms}, we prove Corollary~\ref{cor:BK90_newproof}.
\begin{proof}[Proof of Corollary~\ref{cor:BK90_newproof}]
	For each \(1\le i\le n\), let \(R_i\) be the operator obtained from \(T\) by keeping only the \(i\)-th block row and replacing all other block rows by \(0\). Then
	\[
	T=\sum_{i=1}^n R_i,
	\qquad
	R_i^*R_k=0 \quad (i\ne k),
	\]
	and hence
	\[
	\sum_{i=1}^n |R_i|^2
	=
	\sum_{i=1}^n R_i^*R_i
	=
	T^*T
	=
	|T|^2.
	\]
	Therefore, Theorem~\ref{thm:schatten_p_norms} yields, for \(2\le p<\infty\),
	\begin{equation}\label{eq:row-trunc-corrected}
		n^{1-\frac p2}\|T\|_p^p
		\le
		\sum_{i=1}^n \|R_i\|_p^p
		\le
		\|T\|_p^p.
	\end{equation}
	
	Now fix \(i\), and let
	\[
	\widetilde R_i=
	\begin{pmatrix}
		T_{i1} & T_{i2} & \cdots & T_{in}
	\end{pmatrix}\in \mathbb M_{N,nN}
	\]
	be the corresponding block row operator. Since \(R_i\) is unitarily equivalent to \(\widetilde R_i\oplus 0\), we have
	\[
	\|R_i\|_p=\|\widetilde R_i\|_p.
	\]
	Moreover,
	\[
	|\widetilde R_i^*|^2
	=
	\widetilde R_i\widetilde R_i^*
	=
	\sum_{j=1}^n T_{ij}T_{ij}^*
	=
	\sum_{j=1}^n |T_{ij}^*|^2.
	\]
	Applying Theorem~\ref{thm:schatten_p_norms} to the family
	\(\{T_{ij}^*\}_{j=1}^n\) and the single operator \(\widetilde R_i^*\), we obtain
	\[
	\sum_{j=1}^n \|T_{ij}\|_p^p
	=
	\sum_{j=1}^n \|T_{ij}^*\|_p^p
	\le
	\|\widetilde R_i^*\|_p^p
	=
	\|R_i\|_p^p
	\le
	n^{\frac p2-1}\sum_{j=1}^n \|T_{ij}\|_p^p.
	\]
	Summing over \(i\), we get
	\[
	\sum_{i,j=1}^n \|T_{ij}\|_p^p
	\le
	\sum_{i=1}^n \|R_i\|_p^p
	\le
	n^{\frac p2-1}\sum_{i,j=1}^n \|T_{ij}\|_p^p.
	\]
	Combining this with \eqref{eq:row-trunc-corrected} gives
	\[
	n^{2-p}\|T\|_p^p
	\le
	\sum_{i,j=1}^n \|T_{ij}\|_p^p
	\le
	\|T\|_p^p.
	\]
	
	For \(0<p\le 2\), every application of Theorem~\ref{thm:schatten_p_norms} reverses direction, so both inequalities reverse as well.
\end{proof}
	\section{Proof of Theorem~\ref{thm:schatten_p_q_norms}}\label{sec:proof-spq}
We begin with the classic Hadamard three-lines theorem, see
\cite[p.~133, Problem~3]{SS03}.
	\begin{lem}[\cite{SS03}]\label{lem:three-lines}
		Let $f(z)$ be a bounded function on the strip
		$\{x+iy:\ a\le x\le b\}$
		holomorphic in the interior and continuous on the closed strip. Define
		\[
		M(x):=\sup_{y}\bigl|f(x+iy)\bigr|.
		\]
		Then $\log M(x)$ is convex on $[a,b]$. Equivalently, if $x=ta+(1-t)b$ with $0\le t\le 1$, then
		\[
		M(x)\le M(a)^tM(b)^{1-t}.
		\]
	\end{lem}

If $H\in\mathbb{M}_N$ is positive semidefinite and \(z\in\mathbb C\) with \(\Re z\ge 0\), we define \(H^z\)
by functional calculus, with the convention that \(0^z=0\) on \(\ker H\).
In particular, when \(\Re z=0\), the operator \(H^z\) is unitary on
\(\supp H\) and vanishes on \(\ker H\); hence it is a partial isometry.
This convention will be used throughout the interpolation arguments below. 

The next lemma is the key trace estimate.
	\begin{lem}\label{lem:trace-three-lines}
		Let $1<p\le 2$ and let $q$ be conjugate to $p$.
		Assume $(X_1,\dots,X_m)^T=U(Y_1,\dots,Y_n)^T$ for some matrix $U=(u_{ij})\in \mathbb M_{m,n}$
		satisfying $U^*U\le I_n$ (a contraction).
		Then for any $Z_1,\dots,Z_m\in S_q$,
		\begin{equation}\label{eq:trace-est}
			\left|\tr\!\left(\sum_{i=1}^m Z_iX_i\right)\right|
			\le
			(\max_{i,j}|u_{ij}|)^{\frac{1}{p}-\frac{1}{q}}
			\left(\sum_{j=1}^n \|Y_j\|_p^p\right)^\frac{1}{p}
			\left(\sum_{i=1}^m \|Z_i\|_q^p\right)^\frac{1}{p}.
		\end{equation}
	\end{lem}
	
\begin{proof}
	We may assume finite dimensions; the general case follows by approximation.
	
Write polar decompositions \(Y_j=W_j|Y_j|\) and \(Z_i=V_i|Z_i|\).
For \(z=x+iy\) with \(1/2\le x\le 1\), define
\begin{align}\label{eq:XYZ}
	Y_j(z)&:=
	\begin{cases}
		W_j|Y_j|^{pz}, & Y_j\neq 0,\\
		0, & Y_j=0,
	\end{cases}\nonumber\\
	X_i(z)&:=\sum_{j=1}^n u_{ij}Y_j(z),\\
	Z_i(z)&:=
	\begin{cases}
		\|Z_i\|_q^{\,pz-q(1-z)}V_i|Z_i|^{q(1-z)}, & Z_i\neq 0,\\
		0, & Z_i=0.
	\end{cases}\nonumber
\end{align}
	Note that $Y_j(1/p)=Y_j$ and $Z_i(1/p)=Z_i$.
	Set
	\[
	f(z):=\tr\!\left(\sum_{i=1}^m Z_i(z)X_i(z)\right).
	\]
	Then $f$ is analytic in the interior and continuous on the closed strip.
	Since we are in finite dimensions and $1/2\le \Re z\le 1$, the spectral
	decompositions of $|Y_j|$ and $|Z_i|$ show that the families
	$\{Y_j(z)\}$ and $\{Z_i(z)\}$ are uniformly bounded on the closed strip;
	hence $f$ is bounded there.
		
	\medskip
	\noindent\textbf{Step 1: boundary estimate on \(\Re z=1\).}
Take \(z=1+iy\) with \(y\in\mathbb R\) in \eqref{eq:XYZ}, we have
\begin{align*}
Y_j(1+iy)&=W_j|Y_j|^{p(1+iy)},\\
Z_i(1+iy)
&=
\|Z_i\|_q^{p+i(p+q)y}\,V_i|Z_i|^{-iqy},
\end{align*}
where \(|Y_j|^{ipy}\) is unitary on \(\supp |Y_j|\), the operator
	\(W_j|Y_j|^{ipy}\) is a partial isometry with the same initial space as \(W_j\).
Hence, 
\begin{align*}
		\|Y_j(1+iy)\|_1&=\|Y_j\|_p^p,\\
	\|Z_i(1+iy)\|_\infty&=\|Z_i\|_q^p.
\end{align*}
	Using H\"{o}lder's inequality, we obtain
\begin{equation}\label{eq:trZY-le}
	|\tr(Z_i(1+iy)Y_j(1+iy))|
\le
\|Z_i(1+iy)\|_\infty\,\|Y_j(1+iy)\|_1
=
\|Z_i\|_q^p\,\|Y_j\|_p^p.
\end{equation}
Thus,
	\begin{align}
		|f(1+iy)|
		&=
		\left|\tr\!\left(\sum_{i=1}^m Z_i(1+iy)X_i(1+iy)\right)\right|\nonumber \\
		&=
		\left|\sum_{i=1}^m \tr\!\left(Z_i(1+iy)\sum_{j=1}^n u_{ij}Y_j(1+iy)\right)\right|\quad \text{(by \eqref{eq:XYZ})}\nonumber  \\
		&=
		\left|\sum_{i=1}^m\sum_{j=1}^n u_{ij}\,\tr\!\bigl(Z_i(1+iy)Y_j(1+iy)\bigr)\right| \nonumber \\
		&\le
		\sum_{i=1}^m\sum_{j=1}^n |u_{ij}|\,
		\bigl|\tr\!\bigl(Z_i(1+iy)Y_j(1+iy)\bigr)\bigr| \nonumber \\
		&\le
		\sum_{i=1}^m\sum_{j=1}^n |u_{ij}|\,\|Z_i\|_q^p\,\|Y_j\|_p^p \quad \text{(by \eqref{eq:trZY-le})} \nonumber \\
		&\le
		\Bigl(\max_{i,j}|u_{ij}|\Bigr)
		\left(\sum_{i=1}^m \|Z_i\|_q^p\right)
		\left(\sum_{j=1}^n \|Y_j\|_p^p\right)=:M(1). \label{eq:M-1}
	\end{align}
		\medskip
		\noindent\textbf{Step 2: boundary estimate on $\Re z=1/2$.}
		Take \(z=1/2+iy\) with \(y\in\mathbb R\) in \eqref{eq:XYZ}, then
		\begin{align}\label{eq:YZ-1/2}
			\|Y_j(1/2+iy)\|_2^2&=\|Y_j\|_p^p\nonumber\\
			 \|Z_i(1/2+iy)\|_2^2&=\|Z_i\|_q^p.
		\end{align} 
		By Cauchy--Schwarz inequality in $S_2$,
\begin{equation}\label{Cauchy--Schwarz-inequality-S2}
		|f(z)|
	\le \sum_{i=1}^m \|Z_i(z)\|_2\|X_i(z)\|_2
	\le
	\left(\sum_{i=1}^m\|Z_i(z)\|_2^2\right)^{1/2}
	\left(\sum_{i=1}^m\|X_i(z)\|_2^2\right)^{1/2}.
\end{equation}
Since \(U^*U\le I_n\), the operator \(U\otimes I_{S_2}\) is a contraction on the Hilbert space
\(\ell_2^n(S_2)\). Therefore, from the definitions of \(X_i(z)\) and \(Y_j(z)\) in \eqref{eq:XYZ},
\begin{equation}\label{eq:Xz-le-Yz}
	\sum_{i=1}^m\|X_i(z)\|_2^2
	=
	\bigl\|(U\otimes I_{S_2})(Y_1(z),\ldots,Y_n(z))\bigr\|_{\ell_2^m(S_2)}^2
	\le
	\bigl\|(Y_1(z),\ldots,Y_n(z))\bigr\|_{\ell_2^n(S_2)}^2
	=
	\sum_{j=1}^n\|Y_j(z)\|_2^2.
\end{equation}
		Thus, by \eqref{Cauchy--Schwarz-inequality-S2}, \eqref{eq:Xz-le-Yz} and \eqref{eq:YZ-1/2}, we have
\begin{align}
		|f(1/2+iy)|&\le	\left(\sum_{i=1}^m\|Z_i(1/2+iy)\|_2^2\right)^{1/2}
		\left(\sum_{i=1}^m\|X_i(1/2+iy)\|_2^2\right)^{1/2}\nonumber\\
&\le\left(\sum_{i=1}^m\|Z_i(1/2+iy)\|_2^2\right)^{1/2}
\left(\sum_{j=1}^n\|Y_j(1/2+iy)\|_2^2\right)^{1/2}\nonumber\\
&=\left(\sum_{i=1}^m\|Z_i\|_q^p\right)^{1/2}
	\left(\sum_{j=1}^n\|Y_j\|_p^p\right)^{1/2} \label{eq:M-1/2}
	=:M(1/2).
\end{align}

		\medskip
		\noindent\textbf{Step~3: Interpolate.}
		Apply Lemma~\ref{lem:three-lines} on the strip $[1/2,1]$ at $x=1/p$, \eqref{eq:M-1} and \eqref{eq:M-1/2}:
\begin{align*}
	\left|\tr\!\left(\sum_{i=1}^m Z_iX_i\right)\right|&=|f(1/p)|\\
	&\le M(1)^{2(\frac1p-\frac12)}M(1/2)^{2(1-\frac1p)}\\
	&=	(\max_{i,j}|u_{ij}|)^{\frac{1}{p}-\frac{1}{q}}
	\left(\sum_{j=1}^n \|Y_j\|_p^p\right)^\frac{1}{p}
	\left(\sum_{i=1}^m \|Z_i\|_q^p\right)^\frac{1}{p}
\end{align*}
which is exactly \eqref{eq:trace-est}.
	\end{proof}
	Next, we prove Theorem~\ref{thm:schatten_p_q_norms}.
\begin{proof}[Proof of Theorem~\ref{thm:schatten_p_q_norms}]
	Set
$
	\mu:=\max_{i,j}|u_{ij}|.
$
	Assume first that $1<p\le 2$, and let $q$ be the conjugate exponent of $p$.
For each \(i\), write the polar decomposition \(B_i=U_i|B_i|\) and define
\[
Z_i:=
\begin{cases}
	\|B_i\|_p^{\,q-p}|B_i|^{p-1}U_i^*, & B_i\neq 0,\\
	0, & B_i=0.
\end{cases}
\]
Then \(Z_i\in S_q\), and
\[
\tr(Z_iB_i)=\|B_i\|_p^q,
\qquad
\|Z_i\|_q^p=\|B_i\|_p^q.
\]
 Lemma~\ref{lem:trace-three-lines} with $X_i=B_i$ and $Y_j=A_j$, we get
	\[
	\sum_{i=1}^m\|B_i\|_p^q
	=
	\tr\!\left(\sum_{i=1}^m Z_iB_i\right)
	\le
	\mu^{\frac1p-\frac1q}
	\left(\sum_{j=1}^n \|A_j\|_p^p\right)^\frac{1}{p}
	\left(\sum_{i=1}^m \|B_i\|_p^q\right)^\frac{1}{p}.
	\]
	Hence
	\[
	\sum_{i=1}^m\|B_i\|_p^q
	\le
	\mu^{\frac{q}{p}-1}
	\left(\sum_{j=1}^n \|A_j\|_p^p\right)^{\frac{q}{p}}.
	\]
	
	For the second inequality, note that
	\[
	A=U^*B,
	\]
	and $U^*$ is a contraction with
$
	\max_{j,i}|(U^*)_{ji}|=\mu.
$
	Applying the same argument to the contraction $U^*$ yields
	\[
	\sum_{j=1}^n \|A_j\|_p^q
	\le
	\mu^{\frac{q}{p}-1}
	\left(\sum_{i=1}^m \|B_i\|_p^p\right)^{\frac{q}{p}}.
	\]
	This proves the theorem for $1<p\le 2$.
	
Now let $2\le p<\infty$, and let $q$ be its conjugate exponent. Then $1<q\le 2$.
We prove the reverse of the first inequality; the second one follows in the same way by applying the argument to $A=U^*B$.

By duality of $\ell_p(S_p)$ and $\ell_q(S_q)$,
\[
\left(\sum_{j=1}^n\|A_j\|_p^p\right)^\frac{1}{p}
=
\sup\left\{
\left|\sum_{j=1}^n \tr(W_j^*A_j)\right|:
W_j\in S_q,\ 
\left(\sum_{j=1}^n\|W_j\|_q^q\right)^{\frac{1}{q}}\le 1
\right\}.
\]
Fix such a family $(W_j)_{j=1}^n$, since $A=U^*B$, we have
\[
A_j=\sum_{i=1}^m \overline{u_{ij}}B_i,
\qquad 1\le j\le n.
\]
Hence
\[
\sum_{j=1}^n \tr(W_j^*A_j)
=
\sum_{j=1}^n \tr\!\left(W_j^*\sum_{i=1}^m \overline{u_{ij}}B_i\right)
=
\sum_{i=1}^m \tr\!\left((UW)_i^*B_i\right),
\qquad
(UW)_i:=\sum_{j=1}^n u_{ij}W_j.
\]
By Hölder's inequality for the dual pair $\ell_p(S_q)$--$\ell_q(S_p)$,
\[
\left|\sum_{j=1}^n \tr(W_j^*A_j)\right|
\le
\left(\sum_{i=1}^m\|(UW)_i\|_q^p\right)^\frac{1}{p}
\left(\sum_{i=1}^m\|B_i\|_p^q\right)^{\frac{1}{q}}.
\]
Since $1<q\le 2$, the part already proved (applied with exponent $q$ to the
isometry $U$) gives
\[
\left(\sum_{i=1}^m\|(UW)_i\|_q^p\right)^\frac{1}{p}
\le
\mu^{\frac1q-\frac1p}
\left(\sum_{j=1}^n\|W_j\|_q^q\right)^{\frac{1}{q}}
\le
\mu^{\frac1q-\frac1p}.
\]
Therefore
\[
\left|\sum_{j=1}^n \tr(W_j^*A_j)\right|
\le
\mu^{\frac1q-\frac1p}
\left(\sum_{i=1}^m\|B_i\|_p^q\right)^{\frac{1}{q}}.
\]
Taking the supremum over all such $(W_j)$, we obtain
\[
\left(\sum_{j=1}^n\|A_j\|_p^p\right)^\frac{1}{p}
\le
\mu^{\frac1q-\frac1p}
\left(\sum_{i=1}^m\|B_i\|_p^q\right)^{\frac{1}{q}},
\]
which is equivalent to
\[
\sum_{i=1}^m\|B_i\|_p^q
\ge
\mu^{\frac{q}{p}-1}
\left(\sum_{j=1}^n \|A_j\|_p^p\right)^{\frac{q}{p}}.
\]
\end{proof}
\section{Proof of Theorem~\ref{thm:extension-BCL}}\label{sec:proof-bcl}

Let
\[
\Omega=\{\pm1\}^n
\]
equipped with the uniform probability measure. For a scalar function
\(g\in L_\infty(\Omega)\), write
\[
\mathbb E_\varepsilon[g(\varepsilon)]
=
\frac1{2^n}\sum_{\varepsilon\in\Omega} g(\varepsilon)
\]
for expectation over \(\Omega\).

For \(0\le k\le n\), let
\[
\mathcal F_k=\sigma(\varepsilon_1,\dots,\varepsilon_k)
\]
be the canonical filtration on \(\Omega\), where \(\mathcal F_0=\{\varnothing,\Omega\}\),
and let
\[
\mathbb E_k:L_\infty(\Omega)\to L_\infty(\mathcal F_k)
\]
denote the corresponding scalar conditional expectation.

Set
\[
\mathcal A:=L_\infty(\Omega)\bar\otimes B(\ell_2),
\]
equipped with the product trace
\[
\tau_\Omega\otimes \tr,
\qquad
\tau_\Omega(g)=\mathbb E_\varepsilon[g(\varepsilon)].
\]
For each \(0\le k\le n\), define
\[
\mathcal E_k:=\mathbb E_k\otimes \Id_{B(\ell_2)}.
\]
Then \((\mathcal E_k)_{k=0}^n\) is an increasing family of
trace-preserving conditional expectations on \(\mathcal A\), and
\[
L_p(\mathcal A)\cong L_p(\Omega;S_p).
\]
In particular, if \(X\in L_p(\mathcal A)\) is viewed as an \(S_p\)-valued
function on \(\Omega\), then
\[
\|X\|_{L_p(\mathcal A)}^p
=
\mathbb E_\varepsilon\|X(\varepsilon)\|_{S_p}^p
=
\frac1{2^n}\sum_{\varepsilon\in\{\pm1\}^n}\|X(\varepsilon)\|_p^p.
\]

We shall use the following martingale convexity inequality of Ricard--Xu, see {\cite[Corollary~3]{RX16}}.

\begin{lem}[\cite{RX16}]\label{cor:RX3-finite}
	Let \(1<p\le2\). Then for every \(x\in L_p(\mathcal A)\),
	\[
	\|\mathcal E_n(x)\|_p^2
	\ge
	\|\mathcal E_0(x)\|_p^2
	+
	(p-1)\sum_{k=1}^n \|\mathcal E_k(x)-\mathcal E_{k-1}(x)\|_p^2.
	\]
	For \(2\le p<\infty\), the inequality is reversed.
\end{lem}

\begin{proof}[Proof of Theorem~\ref{thm:extension-BCL}]
	Define
	\[
	X(\varepsilon):=A+\sum_{k=1}^n \varepsilon_k B_k,
	\qquad \varepsilon\in\Omega.
	\]
	Then
	\[
	\|X\|_{L_p(\mathcal A)}^p
	=
	\frac1{2^n}\sum_{\varepsilon\in\{\pm1\}^n}
	\left\|A+\sum_{k=1}^n \varepsilon_k B_k\right\|_p^p.
	\]
	
	Moreover, for each \(0\le k\le n\),
	\[
	\mathcal E_kX(\varepsilon)=A+\sum_{j=1}^k \varepsilon_j B_j.
	\]
	Hence
	\[
	\mathcal E_0X=A,\qquad \mathcal E_nX=X,
	\]
	and
	\[
	\mathcal E_kX-\mathcal E_{k-1}X=\varepsilon_k B_k,\qquad 1\le k\le n.
	\]
	
	Applying Lemma~\ref{cor:RX3-finite} to \(X\), we obtain for \(1<p\le2\),
	\[
	\|X\|_{L_p(\mathcal A)}^2
	=
	\|\mathcal E_nX\|_{L_p(\mathcal A)}^2
	\ge
	\|\mathcal E_0X\|_{L_p(\mathcal A)}^2
	+
	(p-1)\sum_{k=1}^n
	\|\mathcal E_kX-\mathcal E_{k-1}X\|_{L_p(\mathcal A)}^2.
	\]
	For \(2\le p<\infty\), the inequality is reversed.
	
	Now \(\mathcal E_0X\) is the constant function equal to \(A\), so
	\[
	\|\mathcal E_0X\|_{L_p(\mathcal A)}=\|A\|_p.
	\]
	Also, since \(|\varepsilon_k|=1\),
	\[
	\|\mathcal E_kX-\mathcal E_{k-1}X\|_{L_p(\mathcal A)}^p
	=
	\mathbb E_\varepsilon\|\varepsilon_k B_k\|_p^p
	=
	\|B_k\|_p^p,
	\]
	hence
	\[
	\|\mathcal E_kX-\mathcal E_{k-1}X\|_{L_p(\mathcal A)}=\|B_k\|_p.
	\]
	Therefore,
	\[
	\|X\|_{L_p(\mathcal A)}^2
	\ge
	\|A\|_p^2+(p-1)\sum_{k=1}^n \|B_k\|_p^2
	\]
	for \(1<p\le2\), with the reverse inequality for \(2\le p<\infty\).
	
Since the function $t\mapsto t^{p/2}$ is increasing on $[0,\infty)$, raising the
previous inequality to the power $p/2$ yields
\[
\frac1{2^n}\sum_{\varepsilon\in\{\pm1\}^n}
\left\|A+\sum_{k=1}^n \varepsilon_k B_k\right\|_p^p
\ge
\left(\|A\|_p^2+(p-1)\sum_{k=1}^n \|B_k\|_p^2\right)^{p/2}
\]
for $1<p\le2$, and the inequality is reversed for $2\le p<\infty$.
\end{proof}
\section{Proof of Theorem~\ref{thm:Audenaert-extension}}\label{sec:proof-aud}

In this section, we first isolate two special cases in which the optimal constant is \(1\), then record
a duality principle adapted to the present situation, then recall the standard
three-lines interpolation lemma for Schatten $p$-norms, and finally prove
Theorem~\ref{thm:Audenaert-extension}.

We begin with the two endpoint estimates with constant \(1\).

\begin{prop}\label{prop:audext-endpoints}
	Let
$
	T=[A_{ij}]
$ be the $m\times n$ partitioned matrix with $A_{ij}\in \mathbb M_N$ and denote its Schatten \(p\)-norm compression by $C_p(T):=[\|A_{ij}\|_p]$.
	Then
	\begin{enumerate}[label=\textup{(\alph*)}]
		\item
$
		\|C_2(T)\|_2=\|T\|_2.
$
		\item
$
		\|T\|_\infty\le \|C_\infty(T)\|_\infty.
$
	\end{enumerate}
\end{prop}

\begin{proof}
	\textup{(a)} By definition of the Hilbert--Schmidt norm,
	\[
	\|T\|_2^2=\sum_{i=1}^m\sum_{j=1}^n \|A_{ij}\|_2^2.
	\]
	On the other hand, \(C_2(T)\) is the scalar \(m\times n\) matrix whose
	\((i,j)\)-entry is \(\|A_{ij}\|_2\), and therefore
	\[
	\|C_2(T)\|_2^2=\sum_{i=1}^m\sum_{j=1}^n \|A_{ij}\|_2^2.
	\]
	Hence \(\|C_2(T)\|_2=\|T\|_2\).
	
	\medskip
	\noindent\textup{(b)} Let
	\[
	x=(x_1,\dots,x_n)\in (\mathbb C^N)^{\oplus n}\cong \mathbb C^{nN},
	\qquad \|x\|_2=1,
	\]
	where each \(x_j\in \mathbb C^N\). Set
	\[
	a=(\|x_1\|_2,\dots,\|x_n\|_2)^T\in \mathbb R_+^n.
	\]
	Then
	\[
	\|a\|_2^2=\sum_{j=1}^n \|x_j\|_2^2=\|x\|_2^2=1,
	\]
	so \(\|a\|_2=1\).
	
	Now
	\[
	Tx=
	\begin{pmatrix}
		\sum_{j=1}^n A_{1j}x_j\\
		\vdots\\
		\sum_{j=1}^n A_{mj}x_j
	\end{pmatrix},
	\]
	and therefore
	\begin{align*}
		\|Tx\|_2^2
		&=
		\sum_{i=1}^m \left\|\sum_{j=1}^n A_{ij}x_j\right\|_2^2 \\
		&\le
		\sum_{i=1}^m \left(\sum_{j=1}^n \|A_{ij}\|_\infty \|x_j\|_2\right)^2 \\
		&=
		\|C_\infty(T)a\|_2^2 \\
		&\le
		\|C_\infty(T)\|_\infty^2\,\|a\|_2^2
		=
		\|C_\infty(T)\|_\infty^2.
	\end{align*}
	Taking the supremum over all \(x\) with \(\|x\|_2=1\), we obtain
	\[
	\|T\|_\infty\le \|C_\infty(T)\|_\infty.
	\]
\end{proof}

The next proposition is the duality principle that we shall use to pass from
an estimate of the form \(\|C_p(\cdot)\|_2\le K\|\cdot\|_p\) on the range
\(1\le p\le 2\) to the reverse estimate \(\|\cdot\|_q\le K\|C_q(\cdot)\|_2\)
on the conjugate range.

\begin{prop}\label{prop:audext-duality}
	Let \(1\le p\le 2\), and let \(q\) be the conjugate exponent of \(p\).
	Suppose that there exists a constant \(K\ge 0\) such that, for every
	\(m\times n\) partitioned matrix $T=[A_{ij}]$ with $A_{ij}\in \mathbb M_N$ and its Schatten \(p\)-norm compression  $C_p(T):=[\|A_{ij}\|_p]$,
	one has
	\[
	\|C_p(T)\|_2\le K\|T\|_p.
	\]
	Then
	we have
	\[
	\|T\|_q\le K\|C_q(T)\|_2.
	\]
\end{prop}

\begin{proof}
	By Schatten duality on \(\mathbb M_{mN,nN}\),
	\[
	\|T\|_q
	=
	\sup\Bigl\{
	|\Tr(X^*T)|:\ X\in \mathbb M_{mN,nN},\ \|X\|_p\le 1
	\Bigr\}.
	\]
	Fix such an \(X\), and write it in block form as
	\[
	X=[X_{ij}]_{1\le i\le m,\;1\le j\le n},
	\qquad X_{ij}\in \mathbb M_N.
	\]
	Then
	\begin{align*}
		|\Tr(X^*T)|&=\left| \sum_{i=1}^m\sum_{j=1}^n \Tr(X_{ij}^*A_{ij})\right| \\
		&\le
		\sum_{i=1}^m\sum_{j=1}^n \|X_{ij}\|_p\,\|A_{ij}\|_q\quad \text{(By H\"{o}lder's inequality)} \\
		&=
		\left\langle C_p(X),\, C_q(T)\right\rangle_{\mathrm{HS}} \\
		&\le
		\|C_p(X)\|_2\,\|C_q(T)\|_2\quad \text{(By Cauchy-Schwarz inequality)}.
	\end{align*}
	By the assumed estimate,
	\[
	\|C_p(X)\|_2\le K\|X\|_p\le K.
	\]
	Therefore
	\[
	|\Tr(X^*T)|\le K\|C_q(T)\|_2
	\]
	for every \(X\) with \(\|X\|_p\le 1\). Taking the supremum over all such \(X\)
	yields
	\[
	\|T\|_q\le K\|C_q(T)\|_2.
	\]
\end{proof}

To interpolate between the two special cases above, we shall use the following
standard three-lines estimate for Schatten norms; see
\cite[Chapter~4]{BL76} and \cite{Cal64}.

\begin{lem}\label{lem:audext-three-lines}
	Let \(d_1,d_2\in \mathbb N\), and let \(F(z)\) be a function on the strip
	\[
	\{z=x+iy:\ 0\le x\le 1\}
	\]
	which is analytic in the interior and continuous on the closed strip, with
	values in \(\mathbb M_{d_1,d_2}\).
	Define
	\[
	M_0:=\sup_{t\in\mathbb R}\|F(it)\|_\infty,
	\qquad
	M_1:=\sup_{t\in\mathbb R}\|F(1+it)\|_2.
	\]
	Then for every \(\theta\in(0,1)\),
	\[
	\|F(\theta)\|_{2/\theta}\le M_0^{1-\theta}M_1^\theta.
	\]
\end{lem}

We are now ready to prove the core \(2\)-compression estimate, from which
Theorem~\ref{thm:Audenaert-extension} will follow immediately by comparing the
\(\ell_p\)- and  \(\ell_2\)-estimates.

\begin{thm}\label{thm:audext-core-2compression}
	Let
$
T=[A_{ij}]
$ be the $m\times n$ partitioned matrix with $A_{ij}\in \mathbb M_N$ and denote its Schatten \(p\)-norm compression by $C_p(T):=[\|A_{ij}\|_p]$. Then
	\begin{align}
		\|C_p(T)\|_2&\le \|T\|_p, \qquad\qquad  1\le p\le 2,\label{eq:audext-core-12}\\
		\|T\|_p&\le \|C_p(T)\|_2, \qquad 2\le p\le \infty.\label{eq:audext-core-2inf}
	\end{align}
\end{thm}

\begin{proof}
	We first prove \eqref{eq:audext-core-12}.
	
	\medskip
	\noindent\textbf{Step 1: the endpoint \(p=2\).}
	By Proposition~\ref{prop:audext-endpoints}\textup{(a)},
	\[
	\|C_2(T)\|_2=\|T\|_2.
	\]
	So \eqref{eq:audext-core-12} holds at \(p=2\), in fact with equality.
	
	\medskip
	\noindent\textbf{Step 2: the endpoint \(p=1\).}
	If \(C_1(T)=0\), then \(\|A_{ij}\|_1=0\) for all \(i,j\), hence \(A_{ij}=0\)
	for all \(i,j\), so \(T=0\), and there is nothing to prove.
	
	Assume now that \(C_1(T)\neq 0\), and set
	\[
	Y:=\frac{C_1(T)}{\|C_1(T)\|_2}=(y_{ij})\in \mathbb M_{m,n}.
	\]
	Then \(Y\) has nonnegative entries and \(\|Y\|_2=1\).
	
	For each \((i,j)\), write the polar decomposition
	\[
	A_{ij}=U_{ij}|A_{ij}|,
	\]
	and define
	\[
	X_{ij}:=
	\begin{cases}
		y_{ij}U_{ij},& A_{ij}\neq 0,\\[1mm]
		0,& A_{ij}=0.
	\end{cases}
	\]
Let
$
	X:=[X_{ij}] 
$ be the $m\times n$ partitioned matrix.
	Since \(C_\infty(X)=Y\), Proposition~\ref{prop:audext-endpoints}\textup{(b)}
	gives
	\[
	\|X\|_\infty
	\le
	\|C_\infty(X)\|_\infty
	=
	\|Y\|_\infty
	\le
	\|Y\|_2
	=
	1.
	\]
	
	On the other hand,
	\begin{align*}
		\Tr(X^*T)
		&=
		\sum_{i=1}^m\sum_{j=1}^n \Tr(X_{ij}^*A_{ij}) \\
		&=
		\sum_{i=1}^m\sum_{j=1}^n y_{ij}\Tr(|A_{ij}|) \\
		&=
		\sum_{i=1}^m\sum_{j=1}^n y_{ij}\|A_{ij}\|_1 \\
		&=
		\frac{\sum_{i,j}\|A_{ij}\|_1^2}{\|C_1(T)\|_2}
		=
		\|C_1(T)\|_2.
	\end{align*}
	Therefore, by Hölder's inequality for the dual pair \((S_1,S_\infty)\),
	\[
	\|C_1(T)\|_2
	=
	|\Tr(X^*T)|
	\le
	\|X\|_\infty\,\|T\|_1
	\le
	\|T\|_1.
	\]
	So \eqref{eq:audext-core-12} also holds at \(p=1\).
	
	\medskip
	\noindent\textbf{Step 3: the case \(1<p<2\).}
	Let
	\[
	q=\frac{p}{p-1}\in(2,\infty),
	\qquad
	\theta:=\frac{2}{q}=2\Bigl(1-\frac1p\Bigr)\in(0,1).
	\]
	Then
	\[
	\frac{2}{\theta}=q,
	\qquad
	\frac{p\theta}{2}=\frac{p}{q}=p-1.
	\]
	
	Again, if \(C_p(T)=0\), then \(T=0\), so we may assume \(C_p(T)\neq 0\).
	Define
	\[
	Y:=\frac{C_p(T)}{\|C_p(T)\|_2}=(y_{ij})\in \mathbb M_{m,n}.
	\]
	Then \(Y\) has nonnegative entries and \(\|Y\|_2=1\).
	
	For each block \(A_{ij}\), write
	\[
	A_{ij}=U_{ij}|A_{ij}|
	\]
	in polar form. For \(z\) in the strip \(0\le \Re z\le 1\), define
	\[
	X_{ij}(z):=
	\begin{cases}
		y_{ij}\,U_{ij}\dfrac{|A_{ij}|^{\frac{pz}{2}}}{\|A_{ij}\|_p^{\frac{pz}{2}}},
		& A_{ij}\neq 0,\\[3mm]
		0,& A_{ij}=0.
	\end{cases}
	\]
	Let
	\[
	X(z):=[X_{ij}(z)]_{1\le i\le m,\;1\le j\le n}\in \mathbb M_{mN,nN}.
	\]
	Since we are in finite dimensions, \(X(z)\) is analytic in the interior of the
	strip, continuous on the closed strip, and bounded there. Hence
	Lemma~\ref{lem:audext-three-lines} applies to \(X(z)\).
	
	We now compute the trace pairing at \(z=\theta\). For every nonzero block
	\(A_{ij}\),
	\begin{align*}
		\Tr\bigl(X_{ij}(\theta)^*A_{ij}\bigr)
		&=
		y_{ij}\|A_{ij}\|_p^{-\frac{p\theta}{2}}
		\Tr\!\left(|A_{ij}|^{\frac{p\theta}{2}}U_{ij}^*U_{ij}|A_{ij}|\right) \\
		&=
		y_{ij}\|A_{ij}\|_p^{-(p-1)}
		\Tr\!\left(|A_{ij}|^p\right) \\
		&=
		y_{ij}\|A_{ij}\|_p.
	\end{align*}
	For zero blocks the identity is trivial. Summing over all \(i,j\), we obtain
	\begin{equation}\label{eq:audext-core-trace}
		\Tr\bigl(X(\theta)^*T\bigr)
		=
		\sum_{i=1}^m\sum_{j=1}^n y_{ij}\|A_{ij}\|_p
		=
		\frac{\sum_{i,j}\|A_{ij}\|_p^2}{\|C_p(T)\|_2}
		=
		\|C_p(T)\|_2.
	\end{equation}
	
	Next we estimate the two boundary norms.
	
	\smallskip
	\noindent\emph{Left boundary \(\Re z=0\).}
	Let \(z=it\). For each nonzero block \(A_{ij}\), the operator
	\(|A_{ij}|^{\frac{ipt}{2}}\) is unitary on \(\supp|A_{ij}|\), so
	\(U_{ij}|A_{ij}|^{\frac{ipt}{2}}\) is a partial isometry of norm \(1\). Hence
	\[
	\|X_{ij}(it)\|_\infty=y_{ij}.
	\]
	Therefore
	\[
	C_\infty(X(it))=Y.
	\]
	Using Proposition~\ref{prop:audext-endpoints}\textup{(b)}, we get
	\[
	\|X(it)\|_\infty
	\le
	\|C_\infty(X(it))\|_\infty
	=
	\|Y\|_\infty
	\le
	\|Y\|_2
	=
	1.
	\]
	
	\smallskip
	\noindent\emph{Right boundary \(\Re z=1\).}
	Let \(z=1+it\). Since \(|A_{ij}|^{\frac{ipt}{2}}\) is unitary on the support,
	we have
	\[
	\|X_{ij}(1+it)\|_2
	=
	y_{ij}
	\left\|
	U_{ij}\frac{|A_{ij}|^{\frac p2+\frac{ipt}{2}}}{\|A_{ij}\|_p^{\frac p2+\frac{ipt}{2}}}
	\right\|_2
	=
	y_{ij}.
	\]
	Hence
	\[
	C_2(X(1+it))=Y.
	\]
	By Proposition~\ref{prop:audext-endpoints}\textup{(a)},
	\[
	\|X(1+it)\|_2
	=
	\|C_2(X(1+it))\|_2
	=
	\|Y\|_2
	=
	1.
	\]
	
	Applying Lemma~\ref{lem:audext-three-lines} to \(X(z)\), we obtain
	\[
	\|X(\theta)\|_{2/\theta}=\|X(\theta)\|_q\le 1.
	\]
	Combining this with \eqref{eq:audext-core-trace} and Hölder's inequality gives
	\[
	\|C_p(T)\|_2
	=
	|\Tr(X(\theta)^*T)|
	\le
	\|X(\theta)\|_q\,\|T\|_p
	\le
	\|T\|_p.
	\]
	This proves \eqref{eq:audext-core-12} for all \(1\le p\le 2\).
	
	\medskip
	\noindent\textbf{Step 4: proof of \eqref{eq:audext-core-2inf}.}
	Let \(2\le p\le \infty\), and let \(q\in[1,2]\) be the conjugate exponent of
	\(p\). By what we have already proved, every \(m\times n\) partitioned matrix
	\(S\) satisfies
	\[
	\|C_q(S)\|_2\le \|S\|_q.
	\]
	Hence Proposition~\ref{prop:audext-duality} applies with exponent \(q\) and
	constant \(K=1\). Therefore
	\[
	\|T\|_p\le \|C_p(T)\|_2.
	\]
	This is exactly \eqref{eq:audext-core-2inf}.
\end{proof}

We can now deduce Theorem~\ref{thm:Audenaert-extension} immediately.

\begin{proof}[Proof of Theorem~\ref{thm:Audenaert-extension}]
	Set
	\[
	r:=\rank C_p(T)\le \min\{m,n\}.
	\]
	
	If \(1\le p\le 2\), then by Theorem~\ref{thm:audext-core-2compression},
	\[
	\|C_p(T)\|_2\le \|T\|_p.
	\]
	Since \(C_p(T)\) is a scalar matrix of rank \(r\), comparison of Schatten
	norms on its nonzero singular values yields
	\[
	\|C_p(T)\|_p\le r^{\frac1p-\frac12}\|C_p(T)\|_2.
	\]
	Combining the two inequalities, we obtain
	\[
	\|C_p(T)\|_p
	\le
	r^{\frac1p-\frac12}\|T\|_p.
	\]
	
	If \(2\le p\le \infty\), then by
	Theorem~\ref{thm:audext-core-2compression},
	\[
	\|T\|_p\le \|C_p(T)\|_2.
	\]
	Again, since \(C_p(T)\) has rank \(r\),
	\[
	\|C_p(T)\|_2\le r^{\frac12-\frac1p}\|C_p(T)\|_p.
	\]
	Hence
	\[
	\|T\|_p
	\le
	r^{\frac12-\frac1p}\|C_p(T)\|_p.
	\]
	
	This proves Theorem~\ref{thm:Audenaert-extension}.
\end{proof}



\end{document}